\input ./style/arxiv-general.cfg
\documentclass[aap,MSNbibl,secthm,seceqn,dvips]{arximspdf}
\makeatletter
   \@ifpackageloaded{graphicx}{}{\usepackage{graphicx}}
\makeatother

\doi{10.1214/14-AAP1078}
\volume{25}
\issue{6}
\pubyear{2015}
\firstpage{3405}
\lastpage{3433}
\docsubty{FLA}

\makeatletter
\newcommand{\ds}{\displaystyle}
\newcommand{\dd}{d}
\newproclaim{algorithm}{Algorithm}[section]
\newproclaim{rmk}{Remark}[section]
\newproclaim{defi}{Definition}[section]
\newtheorem{prop}[thm]{Proposition}
\newtheorem{lem}[thm]{Lemma}
\makeatother

\begin{document}
\begin{frontmatter}

\title{Robustness of the ${N}$-CUSUM stopping rule in a~Wiener disorder problem}
\runtitle{Robustness of the $N$-CUSUM stopping rule}

\begin{aug}
\author[A]{\fnms{Hongzhong}~\snm{Zhang}\corref{}\ead[label=e1]{hzhang@stat.columbia.edu}},
\author[B]{\fnms{Neofytos}~\snm{Rodosthenous}\ead[label=e2]{N.Rodosthenous@qmul.ac.uk}}
\and
\author[C]{\fnms{Olympia}~\snm{Hadjiliadis}\thanksref{T2}\ead[label=e3]{olympia.hadjiliadis@gmail.com}}
\runauthor{H. Zhang, N. Rodosthenous and O. Hadjiliadis}
\affiliation{Columbia University, Queen Mary University
of London\\ and CUNY}
\address[A]{H. Zhang\\
Department of Statistics\\
Columbia University\\
1255 Amsterdam Ave.\\
New York, New York 10027\\
USA\\
\printead{e1}}
\address[B]{N. Rodosthenous\\
School of Mathematical Sciences \\
Queen Mary University of London\\
London E1 4NS\\
United Kingdom\\
\printead{e2}}
\address[C]{O. Hadjiliadis\\
Department of Mathematics and Statistics\\
Hunter College\\
CUNY\\
919 Hunter East\\
695 Park ave\\
New York, New York 10065\\
USA\\
and\\
Departments of Mathematics\\
\quad and Computer Science\\
Graduate Center\\
CUNY\\
365 5th ave, rm 4208\\
New York, New York 10016\\
USA\\
\printead{e3}}
\end{aug}
\thankstext{T2}{Supported in part by RF-CUNY Collaborative Grant
80209-04 15, PSC-CUNY Grant 65625-00 43, NSA-MSP Grant 081103,
NSF-CCF-MSC Grant 0916452, and NSF-DMS Grant 1222526.}

%
\received{\smonth{10} \syear{2013}}
%
\revised{\smonth{10} \syear{2014}}

%
\begin{abstract}
We study a Wiener disorder problem of detecting the minimum of $N$
change-points in $N$ observation channels coupled by correlated noises.
It is assumed that the observations in each dimension can have
different strengths and that the change-points may differ from channel
to channel. The objective is the quickest detection of the minimum of
the $N$ change-points. We adopt a min--max approach and consider an
extended Lorden's criterion, which is minimized subject to a constraint
on the mean time to the first false alarm. It is seen that, under
partial information of the post-change drifts and a general nonsingular
stochastic correlation structure in the noises, the minimum of~$N$
cumulative sums (CUSUM) stopping rules is asymptotically optimal as the
mean time to the first false alarm increases without bound. We further
discuss applications of this result with emphasis on its implications
to the efficiency of the decentralized versus the centralized systems
of observations which arise in engineering.
\end{abstract}

%
\begin{keyword}[class=AMS]
\kwd[Primary ]{62L10}
\kwd{60K35}
\kwd[; secondary ]{62L15}
\kwd{62C20}
\kwd{60G40}
\end{keyword}
\begin{keyword}
\kwd{CUSUM}
\kwd{correlated noise}
\kwd{quickest detection}
\kwd{Wiener disorder problem}
\end{keyword}
\end{frontmatter}

\section{Introduction}
The problem of quickest detection has been known in the engineering
literature since the 1930s.
Since then there have been various analytical considerations of the
quickest detection problem in
a variety of models and setups (see \cite{PoorHadj} for an overview).
The quickest detection problem, also known as the disorder problem,
concerns the detection
of a change point in the statistical behavior of a stream of sequential
observations. The objective
is to blanace the trade off between a small detection delay and small
frequency of false alarms.
Of this problem there are two main formulations, the Bayesian and the min--max.
In the former, the change point or disorder time is assumed to have an
a priori distribution
usually independent of the observation process while in the latter it
is assumed to be an unknown constant.
An interesting variation of the Bayesian problem in which the change
point is assumed to depend on the observations
is discussed in \cite{Mous08} and treated under Poisson dynamics in
\cite{Seze}.\vadjust{\goodbreak}

Yet in all formulations considered thus far, it is assumed that there
is either one stream of observations
in which there is one \cite{Beib96,Gape,IvanMerz,Mous86,Mous04,Shir96} or multiple alternatives regarding the law of the
post change distribution of the observations \cite{BayrDayaKara,Beib97,BeibLerc}, or alternatively,
multiple streams of observations of various models all undergoing a
disorder at the same time \cite{DayaPoorSeze,Mous06,TartKim06,TartVeer04,TartVeer08}.
In our work, we assume that there are $N$ sources of observations
coupled by correlated noise. The observations are assumed
to be continuous and thus a Wiener model is used. The problem
considered in this work is that in which the $N$ different streams
of observations coupled by correlated noise may undergo a change at $N$
distinct change points.
The objective is then to detect the minimum of the change points or
disorder times. Of this type of problem there has thus far been a
Bayesian formulation in independent streams of Poisson observations
\cite{BayrPoor}. Recently, the case was also considered
of change points that propagate in a sensor array \cite{RaghVeer}.
However, in this configuration the
propagation of the change points depends on the \textit{unknown} identity
of the first sensor affected
and considers a restricted Markovian mechanism of propagation of the change.

In
this paper, we consider the case in which the change points can be
different and do not propagate
in any specific configuration. In fact, in our formulation the change
points or disorder times
are assumed to be unknown constants and a min--max approach to their
estimation is taken.
In particular, we consider an extended Lorden criterion to measure the
worst detection delay over all observation paths and change points.
The objective is then to find a stopping rule that minimizes the
detection delay subject to a lower bound constraint on the mean time to
the first false alarm. The $N$ streams of observations are coupled
through correlated noise. In particular, correlations are modeled
through a
stochastic correlation matrix that is assumed to be nonsingular and
predictable. This work is a continuation of the problem considered in
\cite{HadjZhanPoor} in which the case is considered of independent observations
received at each sensor. In that work, it is seen that the
decentralized system of sensors in which each sensor employs its own
cumulative sum (CUSUM) \cite{PoorHadj} strategy
and communicates its detection through a binary asynchronous message to
a central fusion center, which in turn decides at the first onset of a
signal based on the first communication
performs asymptotically just as well as the centralized system. In
other words, the minimum of $N$-CUSUMs is asymptotically optimal in
detecting the minimum of $N$ distinct
change points in the case of independent observations as the mean time
to the first false alarm increases without bound. The mean time to the
first false alarm can be used as a
benchmark in actual applications in which the engineer or scientist may
make several runs of the system while it is in control in order to
uniquely identify, the appropriate parameter
that would lead to a tolerable rate of false detection. The problem of
optimal detection then boils down to minimizing the detection delay
subject to a tolerable rate of false alarms.
Asymptotic optimality is then proven by comparing the rate of increase
in detection delay to the rate of false alarms as the threshold
parameter varies.
A series of more recent related work includes the case in which the
system of sensors is coupled through the drift parameter as opposed to
the noise \cite{HadjSchaPoor,ZhanHadjSICON}.
In that work, it is once again seen that the minimum of $N$-CUSUMs is
also asymptotically optimal in detecting the minimum of $N$ distinct
change points with respect to a generalized
Kullback--Leibler distance criterion inspired by \cite{Mous04}.

Yet, in none of the above cases is the case of correlated noise
considered even though it is very important in practical applications.
In fact, there are multiple applications of this problem especially in
the area of communications where sensor networks are widely used
and multiple correlated streams of observations are present. The change
points, usually representing the onset of a signal in a specific sensor,
may well be distinct. The minimum of the change points then represents
the onset of a signal in the system. The presence of correlations is due
to the fact that, although sensors are placed typically at different
locations, they are subject to the same physical environment. For example,
in the case of sensors monitoring traffic in opposite (same) directions
may have negative (positive) correlations due to environmental factors
such as the direction
of the wind \cite{CuchPiccMell}. Moreover, the appearance
of a signal at one location may or may not cause interference of the
signal at another location, thus causing
correlations whose structure may even be time or observations
dependent. This happens when the sensors are closely spaced
relative to the curvature of the field being sensed. For example,
temperature sensors or humidity sensors that are in a similar
geographic region will produce readings that are correlated. A
stochastic correlation matrix would best describe such a situation.
Some of the relevant literature that includes such examples can be
found in \cite{BassAbdeBenv,BassBenvGourMeve,BassMeveGour,Ewin,HeylLammSas,Juan,PeetRoec}.

\setcounter{footnote}{1}

In an earlier work \cite{ZhanHadjCorr} the authors treat the
problem of quickest detection of the minimum of two change points in
the special case of two streams of sequential
observations when the correlation in the noise of the observations is
constant and negative and the same drifts are assumed after each of the
disorder times.
This work treats the general case of $N$ correlated streams of
observations in the presence of partial information regarding the
post-change drifts which can as such be different. Moreover, we
consider a general
stochastic correlation matrix allowing for both positive and negative
time and state dependent correlations in the system.
The results found in this work are in fact rather surprising. It is
seen that the minimum of $N$-CUSUM stopping rules maintains its
asymptotically optimal character as the mean time to
the first false alarm increases without bound even in
the case of partially known drifts and a stochastic correlation matrix
coupling the noise of $N$ streams of observations. In particular, it is
proved that the $N$-CUSUM stopping rule (defined in Algorithm~\ref{NCUSUM}) is second-order asymptotically
optimal\footnote{See Definition \ref{optdef} below.} in the case the
post-disorder drift parameters assumed across the $N$ streams of
observations are the same, and is third-order asymptotically optimal
when the post-disorder
drift parameters
are different for an appropriately chosen set of threshold parameters
whose form is explicitly given.

The method used to prove the asymptotic optimality of the $N$-CUSUM
stopping rule is to bound the optimal detection delay from above and
from below. Then we examine the rate at which the difference between
the upper and the lower bounds approach each other as the mean time to
the first false alarm increases without bound. This method is similar
to \cite{GVGMIT,HadjMous,HadjSchaPoor,HadjZhanPoor,Mous04,PT10,ZhanHadjCorr,ZhanHadjSICON}. However, the methodology developed in
this work for establishing the upper and lower bounds is more efficient
and robust in that it is based on probabilistic arguments. In contrast,
the existing work in continuous-time, which is either relied on brute
computation of the asymptotic behaviors of maximum drawdown densities
\cite{HadjZhanPoor} or
on the derivation of sharp solutions to Dirichlet problems with Neumann
conditions \cite{HadjSchaPoor,Mous04,ZhanHadjCorr,ZhanHadjSICON}, is
very difficult in high-dimension and highly sensitive to the
model parameters. The methodology developed in this paper is universal
and can thus
handle a non-Markovian, predictable correlation matrix process for the
noises, which is very useful in practical applications. Finally, our
methodology can be applied to other detection problems not covered in
this paper, for example, quickest detection with multiple alternatives
\cite{O.SQA09,HadjMous}. In establishing the lower bound, we give a
nontrivial generalization of a measure change technique developed in
\cite{Mous04} to $N$-dimensions. Although we do not get the exact
optimality as in one dimension \cite{Mous04}, we do
prove that the optimal detection delay in $N$-dimensions is bounded
from below by that obtained in one dimension, under any predictable,
nonsingular correlation matrix.

In the next section, we formulate the problem mathematically, review
the existing results in one dimension, and introduce the $N$-CUSUM
stopping rule.
In Section~\ref{Robust}, we establish a robust upper bound and a robust
lower bound for both the optimal detection delay and the detection
delay of the $N$-CUSUM stopping rule. These bounds are then used in
Section~\ref{Asympt} to show the main result of the paper---the
asymptotic optimality of the $N$-CUSUM stopping rule under complete or
partial information of the drifts and a stochastic cross-correlated
noise structure in the observations. Applications of these results are
discussed in Section~\ref{applications}.
We conclude with some closing remarks in Section~\ref{conclusion}. The
proof of the lemma that is omitted can be found in the \hyperref[pf]{Appendix}.

Throughout the paper, we denote by $s\wedge t=\min\{s,t\}$, $\mathbb
{R}=(-\infty,\infty)$, $\mathbb{R}_+=[0,\infty)$ and $\overline{\mathbb
{R}}_+=[0,\infty]$.

\section{Formulation of the problem}

Consider a filtered probability space
$(\Omega, \overline{\mathcal{F}}, \overline{\mathbb{F}}, P)$ with
filtration $\overline{\mathbb{F}} = (\overline{\mathcal{F}}_t)_{t\geq0}$,
and the processes $\xi^{(i)}:=\{\xi_t^{(i)}\}_{t \ge0}$, $i=1,\ldots,N$, are assumed to satisfy the following stochastic differential equations:
%
\begin{equation}
\label{Itodynamics} d\xi_t^{(i)} = \mu_i
\mathbf{1}_{\{t\ge\tau_i\}}\,\dd t+\dd w_t^{(i)}.
\end{equation}
Here, $\{\tau_i\}_{1\le i\le N}$ are deterministic but \textit{unknown}
positive constants or $\infty$, $\{\mu_i\}_{1\le i\le N}$ are positive
constants\footnote{The condition can be relaxed. For example, if we
know a priori that $\mu_i<0$ (but not necessarily the value of it),
then we can take $-\xi^{(i)}$ as the $i$th observation process so that
the post-change drift is $-\mu_i>0$. We do not treat in this paper,
however, the case in which we do not know the sign of the post-change
drift.} that are either completely known or partially known. In the
latter case, we assume that $\mu_1>0$ is a known constant, and for
$i=2,\ldots,N$, there are known positive constants\vspace*{0.5pt} $\mu_1\le
\underline{\mu}_i\le\overline{\mu}_i$ such that $\mu_i\in[\underline{\mu
}_i,\overline{\mu}_i]$ holds.\footnote{If $\mu_i$ is known, we can
conveniently take $\underline{\mu}_i=\overline{\mu}_i=\mu_i$.}
The processes $\{w^{(i)}\}_{1\le i\le N}$ for $w^{(i)}:=\{w_t^{(i)}\}
_{t\ge0}$ are
$N$ correlated standard Brownian motions with a predictable,
nonsingular, stochastic instantaneous correlation matrix $\Sigma
_t=(\rho
_t^{i,j})$.
That is, $\rho^{i,j}_{t}$ is the\vspace*{1pt} instantaneous correlation between
Brownian motions $w^{(i)}$ and $w^{(j)}$ (see also \cite{Shre04}, page 227).

An example covered by the above assumptions is one in which $\rho_t^{i,
j}=\rho \mathrm{e}^{-t}$ for $i \ne j$ and some $\rho\in(0,1)$. In other words,
there is a deterministic exponential decay in the instantaneous
correlation of the two sensors $i$ and $j$. Such a situation may arise
by the sudden arrival of a passing rainstorm at sensors $i$ and $j$,
which are customarily placed in the same geographical region and are
therefore also subject to the same climate conditions. Yet our
formulation is even more general in that it is also able to capture
state dependent correlations which is a very realistic scenario since
observations of higher intensity are typically more likely\vspace*{1pt} to cause
higher correlations in the noise, for instance, $N=2$ and $\rho
_t^{i,j}=\frac{\xi_t^{(i)} \xi_t^{(j)}}{1+|\xi_t^{(i)}\xi_t^{(j)}|}
\mathrm{e}^{-t}$ for $i\ne j$.
Another example of a correlated nonstationary white noise structure
arises in the problem of monitoring the vibration of a mechanical
system and is discussed in full detail in Section~11.1.4.1
of \cite{BassNiki}.

To facilitate our analysis, we introduce a family of probability
measures on
the canonical space $(C(\mathbb{R}_+^N),\mathbb{F})$: $\{\mathbb{P}
_{s_1,\ldots,s_N}\}_{(s_1,\ldots,s_N)\in(\overline{\mathbb{R}}_+)^N}$.
Here, $\mathbb{P}_{s_1,\ldots,s_N}$
corresponds to the measure generated on $C(\mathbb{R}_+^N)$ by the processes
$(\xi^{(1)},\ldots,\break \xi^{(N)})$ when the change in the $N$-tuple
process occurs at the time points $\tau_i=s_i$, $1\le i\le N$,
respectively. In particular,
the measure $\mathbb{P}_{\infty,\ldots,\infty}$ characterizes the
law of $N$
correlated standard Brownian motions $\{w^{(i)}\}_{1\le i\le N}$. For
other $s_i$'s, the measure $\mathbb{P}_{s_1,\ldots,s_N}$
can be defined through the Radon--Nikodym derivative process $\frac
{\dd \mathbb{P}_{s_1,\ldots,s_N}}{\dd \mathbb
{P}_{\infty,\ldots,\infty}}
|_{\mathcal{F}_t}$. To this end, we assume that the correlation matrix
$\Sigma_t$ fulfills the Novikov condition:
%
\begin{eqnarray}
\mathbb{E}_{\infty,\ldots,\infty} \biggl\{\exp \biggl(\frac{1}{2}
\biggl\langle\log\biggl(\frac{\dd \mathbb{P}_{s_1,\ldots,s_N}}{\dd \mathbb
{P}_{\infty,\ldots,\infty}} \Big|_{\mathcal{F}_t}\biggr)\biggr\rangle
\biggr) \biggr\}&<& \infty
\nonumber
\\[-8pt]
\label{eq:Nov}
\\[-8pt]
\eqntext{\ds\forall t\ge0, \forall(s_1,
\ldots,s_N)\in(\overline{\mathbb{R}}_+)^N.}
\end{eqnarray}
We comment that the ``reality'' measure $\mathbb{P}_{\tau_1,\ldots
,\tau_N}$ is
one unknown element in $\{\mathbb{P}_{s_1,\ldots,s_N}\}_{(s_1,\ldots
,s_N)\in (\overline{\mathbb{R}}_+)^N}$.

To describe the ``marginal'' law of the $i$th component of the
$N$-tuple process $(\xi^{(1)},\ldots,\xi^{(N)})$, we also introduce
the measure
$\{\mathbb{P}_{s_i}^i\}$, which is the probability measure generated
by the
process $\xi^{(i)}$ on the space
$(C(\mathbb{R}), \mathbb{G}^{(i)})$, where $\mathbb{G}^{(i)}=\{
\mathcal
{G}_t^{(i)}\}_{t \ge0}$ for
$\mathcal{G}_t^{(i)}=\sigma\{(\xi_s^{(i)}); s\le t\}$, is the natural
filtration of $\xi^{(i)}$, and $\tau_i=s_i$ is the value of the
change-point for process $\xi^{(i)}$.

Our objective is to find a stopping rule $T$, which is adapted to the
natural filtration $\mathbb{F}=(\mathcal{F}_t)_{t\ge0}$: $\mathcal
{F}_t=\sigma(\xi_s^{(1)},\ldots,\xi_s^{(N)}; s\le t)$,\footnote{Note
that $\Sigma$ needs not to be adapted to $\mathbb{F}$. For example,
$\Sigma$ can be driven by a $N$-dimensional Brownian motion which is
independent of $w^{(i)}$'s.} to balance the trade-off between
a small detection delay subject to a lower bound on the mean-time to
the first false
alarm and will ultimately detect $\tau_1 \wedge\tau_2 \wedge\cdots
\wedge\tau_N$, which will be denoted by $\tilde{\tau}$ in what follows.
As a performance measure, we consider
%
\begin{equation}
\label{JKL} J^{(N)}(T) =\mathop{\sup_{(s_1, \ldots, s_N ) \in\overline{\mathbb{R}}_+^N}}_{
\tilde{s}<\infty}
\operatorname{essup}\mathbb{E}_{s_1,\ldots, s_N} \bigl\{ (T-\tilde{s})^+|\mathcal{F}_{\tilde{s}} \bigr\},
\end{equation}
where $\tilde{s}=s_1\wedge s_2\wedge\cdots\wedge s_N$, $\mathbb
{E}_{s_1,\ldots, s_N}$ denotes the expectation
under the probability measure $\mathbb{P}_{s_1,\ldots, s_N}$, and the supremum
over $s_1, \ldots, s_N$ is taken over the set in which
$\tilde{s}<\infty$. In other words, we consider the worst detection
delay over all possible realizations of
paths of the $N$-tuple of the stochastic process $\{(\xi_t^{(1)},
\ldots,
\xi_t^{(N)}) \}_{t\ge0}$ up to time $\tilde{s}$, and then consider
the worst
detection delay over all possible
$N$-tuples 
$(s_1, \ldots, s_N)$ over a
set in which at least one of the components takes a finite value.
This is because $T$ is a stopping rule meant to detect the minimum of
the $N$
change points and, therefore, if one of the $N$ processes undergoes a regime
change, any unit of time by which $T$ delays in reacting, should be counted
toward the detection delay. Although it seems to be a quite pessimistic
measure for detection delay, this framework has the merit that one does
not need to impose any prior knowledge of the distribution of the
change-points $\tau_i$'s, as discussed in \cite{Mous08}.
In all, this gives rise to the following stochastic optimization problem:
%
\begin{equation}
\label{eqnproblem}
\quad\inf_{T\in\mathcal{T}_\gamma}J^{(N)}(T)
\qquad \mbox{with  }\mathcal{T}_\gamma= \bigl\{\mathbb{F}\mbox{-stopping rule
}T\dvtx \mathbb{E}_{\infty, \ldots,
\infty}\{ T\} \geq \gamma\bigr\},
\end{equation}
where
$\mathbb{E}_{\infty, \ldots, \infty}\{T\}$ captures the mean time
to the first false
alarm and as such the above constraint describes the tolerance on the
false alarms. In particular, the constant $\gamma>0$ is the lowest
acceptable value of the mean time to the first false alarm. In other
words, the reciprocal of $\gamma$, namely $\frac{1}{\gamma}$, captures
the highest tolerance to the frequency of false alarms of the family of
stopping times considered in this problem.

When detecting $\tau_i$ is our only concern, and that $\mu_i$ is a
known constant,
the problem reduces to an one-dimensional problem of detecting a
one-sided change in a sequence of Brownian motion observations, whose
optimality was found in
\cite{Beib96} and \cite{Shir96}.
It is shown that the optimal stopping rule under Lorden's criterion is
the continuous time version of Page's CUSUM stopping rule, namely the
first passage time of the process
\begin{eqnarray}
\label{1dimpr}
\tilde{y}_t^{(i)}
&:=&  \sup_{0\le s \le t}\dfrac{\dd \mathbb{P}_{s}^i}{{d}\mathbb{P}_\infty^i} \bigg|_{\mathcal{G}_t^{(i)}} =
\tilde{u}_t^{(i)} - \tilde{m}_t^{(i)}
\nonumber
\\[-8pt]
\\[-8pt]
\eqntext{\ds\mbox{for }\tilde{u}_t^{(i)} :=
\mu_i \xi_t^{(i)} - \frac{1}{2}
\mu_i^2 t \quad\mbox{and}\quad \tilde{m}_t^{(i)}
:= \inf_{0\le s\le t} \tilde{u}_s^{(i)},}
\end{eqnarray}
and the CUSUM stopping rule with the threshold $\nu^\star_i >0$ is
given by
%
\begin{equation}
\label{1dimCUSUMstop}
\tilde{T}_{\nu^\star_i}^i = \inf\bigl\{t\ge0 \dvtx
\tilde{y}_t^{(i)} \ge \nu^\star_i\bigr\}.
\end{equation}
The optimal threshold $\nu^\star_i$ is chosen so that
%
\begin{equation}
\label{g}
\hspace*{6pt}\mathbb{E}_\infty^i \bigl\{
\tilde{T}_{\nu^\star_i}^i \bigr\} = \bigl( 2 /
\mu_i^2 \bigr) g\bigl(\nu^\star_i
\bigr) = \gamma\qquad  \mbox{where }g(\nu) := \mathrm{e}^\nu- \nu- 1, \forall
\nu>0.
\end{equation}
The corresponding optimal detection delay achieved by the CUSUM
stopping rule
$\tilde{T}_{\nu^\star_i}^i$ is then given by
%
\begin{equation}
\label{1dimsol}
J^{(1)}\bigl(\tilde{T}_{\nu^\star_i}^i\bigr)
= \mathbb{E}_0^i \bigl\{ \tilde {T}_{\nu^\star
_i}^i
\bigr\} = \frac{2}{\mu_i^2} g\bigl(- \nu^\star_i\bigr).
\end{equation}
The fact that the worst detection delay in the one-dimensional problem
is the same as that incurred in the case
that the change point is exactly at $0$ is a consequence of the nonnegativity
and strong Markov property of the CUSUM process, from which it follows
that the
worst detection delay occurs when the CUSUM process is at $0$ at the
time of the
change (see also \cite{Mous86}).

The optimality of the CUSUM stopping rule in the presence of only one
observation process with a known drift suggests that a CUSUM type of
stopping rule might display
similar optimality properties in the case of multiple observation
processes for
the problem (\ref{eqnproblem}).
In particular, an intuitively appealing rule, when the detection of
$\tilde{\tau}$ is of interest, is to take the minimum of $N$-CUSUM-like
stopping rules (see, e.g., \cite{O.SQA09}), which we formalize in the
following algorithm.

%
\begin{algorithm}\label{NCUSUM}
The $N$-CUSUM\vspace*{1pt} stopping rule with a threshold vector $\hbar=(h_1,\ldots
,h_N)\in(\mathbb{R}_+)^N$ is given by $T_\hbar=T_{h_1}^1\wedge
T_{h_2}^2\wedge\cdots\wedge T_{h_N}^N$, where for each $i=1,\ldots,N$,
%
\begin{eqnarray}
T_{h_i}^{i} &=&
\inf\bigl\{t \ge0 \dvtx  y_t^{(i)} \geq h_i \bigr\}
\nonumber
\\[8pt]
\label{CUSUM1chart}
\\[-25pt]
\eqntext{\ds\mbox{with } y_t^{(i)}=u_t^{(i)}-m_t^{(i)}, \mbox{for } u_t^{(i)}:=\underline {\mu }_i
\xi_t^{(i)}-\frac{1}{2}\underline{
\mu}_i^2t}\\
\eqntext{\ds\mbox{and }m_t^{(i)}:=\inf
_{0\le s\le t}u_t^{(i)}.}
\end{eqnarray}
\end{algorithm}

That is, we use what is known as a multichart CUSUM stopping rule
\cite{Mous06}, which can be written as
\[
T_{\hbar}  =  \inf \biggl\{t \ge0\dvtx  \max \biggl\{
\frac{y_t^{(1)}}{h_1},\ldots,\frac{y_t^{(N)}}{h_N} \biggr\} \ge1 \biggr\},
\]
where $\{y_t^{(i)}\}_{t\ge0}$ is the semi-martingale defined in (\ref{CUSUM1chart}),
for $i=1,\ldots,N$.
We notice that each of the $T_{h_i}^{i}$, for $i=1,\ldots,N$, are
stopping rules also with respect to
each of the smaller filtrations $\mathbb{G}^{(i)}$, and thus they can
be employed by each one of the sensors $S_i$, for each $i$ independently.
Each of the sensors can then subsequently communicate an alarm to a
central fusion
center once its threshold, say $h_i$, is reached by its own CUSUM
statistic process
$y^{(i)}$.
The resulting rule, namely Algorithm \ref{NCUSUM}, can then be devised
by the
central fusion center in that it will declare a detection at the first instance
one of the $N$ sensors communicates.

%
\begin{rmk}
From (\ref{1dimpr}) and (\ref{CUSUM1chart}), it is easily seen that
$y^{(i)}\equiv\tilde{y}^{(i)}$ and $T_{h_i}^i=\tilde{T}_{h_i}^i$,
a.s., provided that $\mu_i=\underline{\mu}_i$ is known. In particular,
we always have $y^{(1)}\equiv\tilde{y}^{(1)}$ and $T_{h_1}^1=\tilde
{T}_{h_1}^1$.
\end{rmk}

While it seems prohibitively difficult to devise a stopping rule that
achieves the optimal detection delay $\inf_{T\in\mathcal{T}_\gamma
}J^{(N)}(T)$ under a general nonsingular correlation matrix $(\Sigma
_t)_{t\ge0}$,
the above $N$-CUSUM stopping rule $T_\hbar$ provides a low-complexity
candidate detection rule for detecting $\tilde{\tau}$.

In particular, we will show that the $N$-CUSUM stopping rule is
asymptotically optimal. To this effect, we give the following
definitions of asymptotic optimality as in \cite{GVGMIT}.
%

\begin{defi}\label{optdef}
Given $\gamma>0$ and a stopping time
$T'\in
\mathcal{T}_\gamma$, we say that:
\begin{enumerate}[3.]
\item[1.] $T'$ has the first-order asymptotic optimality for problem
(\ref{eqnproblem}) if and only if
\[
\lim_{\gamma\to\infty}\frac{J^{(N)}(T')}{\inf_{T\in\mathcal
{T}_\gamma
}J^{(N)}(T)}=1\quad \mbox{and}\quad \lim
_{\gamma\to\infty}\inf_{T\in
\mathcal
{T}_\gamma}J^{(N)}(T)=\infty.
\]
\item[2.] $T'$ has the second-order asymptotic optimality for problem
(\ref{eqnproblem}) if and only if
\[
\lim_{\gamma\to\infty}\Bigl[J^{(N)}\bigl(T'\bigr)-
\inf_{T\in\mathcal{T}_\gamma
}J^{(N)}(T)\Bigr]<\infty\quad\mbox{and}\quad \lim
_{\gamma\to\infty}\inf_{T\in
\mathcal
{T}_\gamma}J^{(N)}(T)=\infty.
\]
\item[3.] $T'$ has the
third-order asymptotic optimality for problem (\ref{eqnproblem}) if and only if
\[
\lim_{\gamma\to\infty}\Bigl[J^{(N)}\bigl(T'\bigr)-
\inf_{T\in\mathcal{T}_\gamma
}J^{(N)}(T)\Bigr]=0\quad \mbox{and}\quad\lim
_{\gamma\to\infty}\inf_{T\in
\mathcal
{T}_\gamma}J^{(N)}(T)=\infty.
\]
\end{enumerate}
\end{defi}

Below we will investigate the performance of $T_\hbar$ by contrasting
it with the optimal detection delay.

\section{Robust bounds for the optimal detection delay} \label{Robust}

In this section, we examine the performance of the $N$-CUSUM stopping
rule by
presenting an upper bound and a lower bound for both the detection
delay of the $N$-CUSUM stopping rule and the
optimal detection
delay defined in (\ref{eqnproblem}). 
To this end, we derive a robust upper bound for the detection delay of
a particular $N$-CUSUM stopping rule~$T_{\hbar}$ in $\mathcal{T}_\gamma$.
Because $T_\hbar$ cannot beat the unknown optimal stopping rule (if
it ever exists), this upper bound will also bound the optimal detection
delay from above.
We then demonstrate that the optimal detection delay in the
$N$-dimensional system is bounded from below by the optimal delay in
1-dimensional systems.


\subsection{The upper bound}
In this subsection, we derive a robust upper bound for the detection
delay of a $N$-CUSUM stopping rule $T_\hbar$, whose thresholds set
$\hbar$ is chosen so that $T_\hbar\in\mathcal{T}_\gamma$ for any
$\gamma>0$.
The upper bound, that we obtain, also dominates the optimal detection
delay, due to the fact that $J^{(N)}(T_\hbar)\ge\inf_{T\in\mathcal
{T}_\gamma}J^{(N)}(T)$ holds.

Now let us introduce
%
\begin{equation}
\label{JjN}
J_j^{(N)}(T)  =  \mathop{\sup_{(s_1, \ldots, s_N) \in\overline{\mathbb{R}}_+^N}}_{
s_j = \tilde{s} < \infty} \operatorname{essup} \mathbb{E}_{s_1,\ldots,s_N}
\bigl\{ (T - s_j)^+ | \mathcal{F}_{s_j} \bigr\},
\end{equation}
for $j=1,\ldots,N$, where $J_j^{(N)}(T)$ is the\vspace*{1pt} detection delay of the
stopping rule $T$ when $s_j \le\min_{i\neq j} \{s_i\}$, implying that the
performance measure defined in (\ref{JKL}) is given by
$J^{(N)}(T) = \max_{1\leq j\leq N} J_j^{(N)}(T)$.
We now consider the case when all drifts $\mu_i$'s are known constants.
In this case, we select $\hbar$ such that
%
\begin{equation}
\label{hi}
\mathbb{E}_0^1\bigl\{T_{h_1}^1
\bigr\}=\mathbb{E}_0^2\bigl\{T_{h_2}^2
\bigr\}=\cdots =\mathbb{E}_0^N\bigl\{T_{h_N}^N
\bigr\},
\end{equation}
or equivalently [by (\ref{1dimsol})],
\[
\frac{1}{\mu_1^2}g(-h_1)=\frac{1}{\mu_2^2}g(-h_2)=
\cdots=\frac
{1}{\mu
_N^2}g(-h_N).
\]
Due to the monotonicity of function $g$, $h_i$'s are uniquely
determined once $h_1>0$ is given. In general, if we only have partial
information about $\mu_i$'s for $i=2,\ldots,N$, we instead consider
%
\begin{equation}
\label{hi1}
\frac{1}{\mu_1^2}g(-h_1)=\frac{1}{\underline{\mu}_2^2}g(-h_2)=
\cdots =\frac{1}{\underline{\mu}_N^2}g(-h_N).
\end{equation}
By choosing the $N$-CUSUM stopping rule $T_{\hbar}$ in this way, we are
able to
get an easily computable upper bound for the worst detection delay
$J^{(N)}(T_\hbar)$. The assertion is proved in the following proposition.
%




%
\begin{prop} \label{UpB}
Suppose that $\hbar\in\mathbb{R}_+^N$ satisfies the equations in
(\ref{hi})
or~(\ref{hi1}), then we have for the $N$-CUSUM stopping rule $T_{\hbar
}$, that
%
\begin{equation}
\label{UpBeq}
J^{(N)}(T_\hbar) \leq\mathbb{E}_0^1
\bigl\{T_{h_1}^1\bigr\}=\frac{2}{\mu_1^2}
g(-h_1),
\end{equation}
where the function $g$ is defined in (\ref{g}).
\end{prop}

\begin{pf}
For any $(s_1, \ldots, s_N) \in\overline{\mathbb{R}}_+^N$ such that $s_j
\le\min_{i\neq j} \{s_i\}$ and $s_j<\infty$, we have
%
\begin{eqnarray}
\mathbb{E}_{s_1,\ldots,s_N} \bigl\{ (T_\hbar-
s_j)^+ | \mathcal {F}_{s_j} \bigr\} &=&  \mathbb{E}_{s_1,\ldots,s_N}
\bigl\{ \bigl(T^{1}_{h_1} \wedge\cdots\wedge
T^{N}_{h_N} - s_j\bigr)^+ | \mathcal{F}_{s_j}
\bigr\}
\nonumber
\\
\label{E<}
&\leq & \mathbb{E}_{s_1,\ldots,s_N} \bigl\{ \bigl(T^{j}_{h_j} -
s_j\bigr)^+ | \mathcal {F}_{s_j} \bigr\} \\
&=&
\mathbb{E}^{j}_{s_j} \bigl\{ \bigl(T^{j}_{h_j}
- s_j\bigr)^+ | \mathcal {G}_{s_j}^j \bigr\},
\nonumber
\end{eqnarray}
where the last equality follows from the fact that the CUSUM stopping rule
$T^{j}_{h_j}$ is $\mathbb{G}^{(j)}$-measurable
and we can thus use the ``marginal'' law of the $j$th component of the $N$-tuple
process $(\xi^{(1)},\ldots,\xi^{(N)})$, given in this case by the measure
$\mathbb{P}_{s_j}^{j}$.
By taking the essential supremum and then the supremum over $s_1,
\ldots, s_N$
such that $s_j \le\min_{i\neq j} \{s_i\}$ on both sides of (\ref{E<}),
and using the definitions in (\ref{JKL}) and (\ref{JjN}) for $N=1$, we
get that
%
\begin{equation}
\label{J1<}
J_j^{(N)}(T_\hbar) \leq
J^{(1)}\bigl(T^{j}_{h_j}\bigr) = \sup
_{s_j<\infty
}\operatorname{essup}\mathbb{E}^{j}_{s_j}
\bigl\{ \bigl(T^{j}_{h_j}-s_j\bigr)^+|\mathcal
{G}_{s_j}^j \bigr\}.
\end{equation}
To get the conditional expectation in the above expression, we use the
strong Markov property of the\vspace*{1pt} processes $y^{(j)}_t$
(see, e.g., \cite{Okse}, Theorem 7.2.4)
and
apply It\^o's formula to $\{g(-y_t^{(j)})\}_{s_j \le t < T_{h_j}^j}$
(see, e.g., \cite{Okse}, Theorem 4.1.2)
for the function $g$ given in (\ref{g})
(see also Shiryaev \cite{Shir96} and Moustakides \cite{Mous04}),
we obtain that (by the monotonicity of $g$)
%
\begin{eqnarray}
g(-h_j) &\ge & \mathbf{1}_{\{T_{h_j}^j>s_j\}
}\bigl[g\bigl(-y_{T_{h_j}^j}^{(j)}
\bigr)-g\bigl(-y_{s_j}^{(j)}\bigr)\bigr]
\nonumber
\\
\label{sing}
&=& \mathbf{1}_{\{T_{h_j}^j>s_j\}} \biggl[\int_{s_j}^{T_{h_j}^j}
\underline{\mu}_j \biggl(\mu_j-\frac{1}{2}
\underline{\mu}_j \biggr)\,\dd s\\
&&\nonumber\hspace*{43pt}{}-\int_{s_j}^{T_{h_j}^j}g'
\bigl(-y_s^{(j)}\bigr)\,\dd m_s^{(j)}+M_{T_{h_j}^j}-M_{s_j}
\biggr],
\end{eqnarray}
where the process $m^{(j)}$ is given by (\ref{CUSUM1chart}) and the
continuous square integrable martingale $M=\{M_{t}\}_{t\ge0}$ is given by
\[
M_{t}=\int_0^{t\wedge T_{h_j}^j}\underline{\mu
}_jg'\bigl(-y_s^{(j)}\bigr)\,\dd
w_s^{(j)}.
\]
Taking into account that the process $m^{(j)}$ decreases only on
the random set  $\{ t \geq0 \dvtx  y_{t}^{(j)} =0 \}$ and the measure
$\dd m_t^{(j)}=0$ off this set, together with the fact that
$g'(0) = 0$,
we conclude that the integral in (\ref{sing}) can be set
equal to zero. We then take the conditional expectations with respect
to the probability measure $\mathbb{P}_{s_j}^j$ given\vspace*{1pt} $\mathcal{G}_{s_j}^j$ in
(\ref{sing}) and by means of the Doob's optional sampling
theorem
(see, e.g., \cite{KS}, Chapter~1, Theorem~3.22), we have
%
\begin{eqnarray}
g(-h_j) &\ge & \mathbf{1}_{\{T_{h_j}^j>s_j\}}\mathbb{E}_{s_j}^j
\biggl\{ \int_{s_j}^{T_{h_j}^j}\underline{
\mu}_j \biggl(\mu_j-\frac
{1}{2}\underline{\mu
}_j \biggr)\,\dd s \Big| \mathcal{G}_{s_j}^j \biggr\}
\nonumber
\\[-8pt]
\label{J11}
\\[-8pt]
\nonumber
&\ge &\frac{\underline{\mu
}_j^2}{2}\mathbb{E}_{s_j}^{j}\bigl\{
\bigl(T_{h_j}^j-s_j\bigr)^+|
\mathcal{G}_{s_j}^j\bigr\}.
\end{eqnarray}
Therefore, by (\ref{J1<}) and (\ref{J11}) we have $J_j^{(N)}(T_\hbar)\le
\frac{2}{\underline{\mu}_j^2}g(-h_j)$.
By the arbitrariness of $j$, we have
%
\begin{eqnarray}
J^{(N)}(T_\hbar) &=&  \max_{1\leq j\leq N}
J_j^{(N)}(T_\hbar) 
\nonumber
\\[-8pt]
\label{JN}\\[-8pt]
\nonumber
&\le&  \max \biggl\{
\frac{2}{\mu_1^2} g(-h_1),\ldots,\frac{2}{\underline
{\mu
}_N^2}
g(-h_N) \biggr\} = \frac{2}{\mu_1^2}g(-h_1),
\end{eqnarray}
where the last equality is a consequence of the equations in (\ref{hi1}).
\end{pf}

Condition (\ref{hi1}) reduces the thresholds' selection problem
from $N$ dimension to one dimension.
In order to bound the optimal detection delay in (\ref{eqnproblem})
using the result in Proposition \ref{UpB}, we will choose $h_1$ so that
the resulting $N$-CUSUM stopping rule $T_\hbar\in\mathcal{T}_\gamma$.
That is,
%
\begin{equation}
\mathbb{E}_{\infty,\ldots,\infty}\{T_\hbar\}\ge\gamma.
\end{equation}
To this end, we derive a lower bound for the mean time to the first
false alarm $\mathbb{E}_{\infty,\ldots,\infty}\{T_\hbar\}$, which
is robust
with respect to the covariance matrix $(\Sigma_t)_{t\ge0}$. In the
sequel, we first study the case of equal drift size with complete
information, where $0 < \mu_1 = \mu_2 = \cdots= \mu_N = \mu$ are known
constant, and then treat the case of unequal drift size with complete
information, such that $\mu_i$'s are all known and
$0 < \mu_1 = \mu_2 = \cdots= \mu_k< \min_{i>k} \mu_i$ holds for some
$k\in\{1,\ldots,N-1\}$. Finally, we study the general case with
partial information, where we only know $\mu_1$ and
the intervals
$[\underline{\mu}_i,\overline{\mu}_i]\ni\mu_i$ for all
$i=2,\ldots,N$.

\subsubsection{Equal drift case---complete information about \texorpdfstring{$\mu_i$}{mui}'s}
In this case, it is assumed that all $\mu_i$ are known and $\mu
_i=\mu
_1=\mu>0$ for all $i=1,\ldots,N$.
Then the monotonicity of function $g$ and (\ref{hi}) imply that $h_1 =
h_2 = \cdots
= h_N = h$ for the $N$-CUSUM stopping rule. Hence, with a slight abuse
of notation, we denote by $T_h=T_{\hbar}$. Below we derive a lower
bound for the
mean time of the first false alarm of the latter. 

%
\begin{prop}\label{thm:false_alarm}
Suppose that all thresholds of the $N$-CUSUM stopping rule are chosen
to be equal to $h>0$.
Then the first false alarm for the $N$-CUSUM stopping rule $T_h$ satisfies
%
\begin{equation}
\label{LBsame}
\mathbb{E}_{\infty,\ldots,\infty}\{T_h\}\ge\frac{1}{N}
\mathbb {E}_\infty^1\bigl\{T_h^1\bigr
\} =\frac{2}{N\mu^2} g(h),
\end{equation}
where the function $g$ is defined in (\ref{g}).
\end{prop}

\begin{pf}
For any $i=1,\ldots,N$, we have
\begin{eqnarray}
\qquad\quad\mathbb{E}_\infty^i \bigl\{
T_h^{i} \bigr\} &=& \mathbb{E}_{\infty,\ldots
,\infty} \bigl\{
T_h^{i} \bigr\} = \mathbb{E}_{\infty,\ldots,\infty}\{
T_h \} + \mathbb{E}_{\infty,\ldots,\infty} \bigl\{ \bigl(T_h^{i}
- T_h\bigr) \mathbf {1}_{\{T_h^{i} \neq T_h\}} \bigr\}
\nonumber
\\[-8pt]
\label{R1}
\\[-8pt]
\nonumber
& = & \mathbb{E}_{\infty,\ldots,\infty}\{ T_h \} + \mathbb{E}_{\infty,\ldots,\infty}
\bigl\{ \mathbb{E}_{\infty
,\ldots,\infty}\bigl\{ T_h^{i} -
T_h | \mathcal{F}_{T_h} \bigr\} \mathbf{1}_{\{T_h^{i} \neq T_h\}}
\bigr\},
\nonumber
\end{eqnarray}
where the third equality follows from the tower property of the conditional
expectation and the finiteness of $T_\hbar$.

As in the proof of Proposition \ref{UpB}, we
apply It\^o's formula to $\{g(y_t^{(i)})\}_{T_h \le t < T_h^i}$
(see, e.g., \cite{Okse}, Theorem~4.1.2)
to obtain that
%
\begin{equation}
\label{Ito's expansion}
\qquad g \bigl(y^{(i)}_{T_h^{i}} \bigr) - g
\bigl(y^{(i)}_{T_h} \bigr) = \frac
{\mu
^2}{2}
\bigl(T_h - T_h^{i}\bigr) 
- \int
_{T_h}^{T_h^{i}} g'\bigl(y^{(i)}_s
\bigr) \,\dd m_s^{(i)} + M_{T_h^{i}} - M_{T_h},
\end{equation}
where the process $m^{(i)}$ is given by (\ref{1dimpr}) and the
continuous square
integrable martingale $M = \{M_{t}\}_{t\geq0}$ (with respect to
$\mathbb{P}_{\infty,\ldots,\infty}$) is given by
%
\begin{equation}
\label{mg}
M_{t} = \mu\int_{0}^{t\wedge T_{h}^i}
g'\bigl(y^{(i)}_s\bigr) \,\dd
w_s^{(i)}.
\end{equation}
Taking into account that the process $m^{(i)}$ decreases only on
the random set \mbox{$\{ t \geq 0 \dvtx   y_{t}^{(i)} =0 \}$} and the measure
$\dd m_t^{(i)}=0$ off this set, together with the fact that
$g'(0) = 0$,
we conclude that the integral in (\ref{Ito's expansion}) can be set
equal to zero. We then take the conditional expectations with respect
to the probability measure $\mathbb{P}_{\infty,\ldots,\infty}$ in
(\ref{Ito's expansion}) and by means of the Doob's optional sampling
theorem
(see, e.g., \cite{KS}, Chapter~1, Theorem~3.22),
we have
%
\begin{equation}
\label{Ito}
\mathbb{E}_{\infty,\ldots,\infty} \bigl\{ g \bigl(y^{(i)}_{T_h^{i}}
\bigr) - g \bigl(y^{(i)}_{T_h} \bigr) | \mathcal
{F}_{T_h} \bigr\} = \frac{\mu^2}{2}\mathbb{E}_{\infty,\ldots,\infty} \bigl\{
T_h^{i} - T_h | \mathcal{F}_{T_h}
\bigr\}.
\end{equation}
%
Therefore, using equation (\ref{Ito}) in the expression of (\ref{R1})
we have that
\begin{eqnarray}
\frac{2}{\mu^2}g(h)&= & \mathbb{E}_\infty^i
\bigl\{ T_h^{i} \bigr\}\nonumber
\\
&=&  \mathbb{E}_{\infty,\ldots,\infty}\{ T_h \} + \mathbb{E}_{\infty,\ldots,\infty}
\biggl\{ \frac{2}{\mu^2} \bigl( g\bigl(y^{(i)}_{T_h^{i}}\bigr) -
g\bigl(y^{(i)}_{T_h}\bigr)\bigr) \mathbf{1}_{\{T_h^{i}
\neq T_h\}}
\biggr\}
\nonumber
\\[-8pt]
\label{R3}
\\[-8pt]
\nonumber
&=&  \mathbb{E}_{\infty,\ldots,\infty}\{ T_h \} + \mathbb{E}_{\infty,\ldots,\infty}
\biggl\{ \frac{2}{\mu^2} \bigl(g(h) - g\bigl(y^{(i)}_{T_h}
\bigr)\bigr) \mathbf{1}_{\{T_h^{i} \neq T_h\}} \biggr\}
\nonumber
\\
&\leq & \mathbb{E}_{\infty,\ldots,\infty}\{ T_h \} + \frac{2}{\mu^2} g(h)
\mathbb{P}_{\infty,\ldots,\infty} \bigl( T_h^{i} \neq
T_h \bigr),
\nonumber
\end{eqnarray}
where the first equality and the function $g$ are given by (\ref{g})
and the
third equality follows
from the definition of the one-dimensional CUSUM stopping rule in
(\ref{CUSUM1chart}).
It follows that
\[
\frac{\mu^2}{2} \mathbb{E}_{\infty,\ldots,\infty}\{T_h\} \ge g(h)
\mathbb{P}_{\infty,\ldots,\infty}\bigl(T_h=T_h^{i}
\bigr).
\]
Hence, by summing both sides over all $i=1,\ldots,N$, we get
\begin{eqnarray*}
\frac{N\mu^2}{2} \mathbb{E}_{\infty,\ldots,\infty
}\{T_h\} &\ge&  g(h) \sum
_{i=1}^N \mathbb{P}_{\infty,\ldots,\infty}
\bigl(T_h = T_h^{i}\bigr)
\\
&\ge & g(h)\mathbb{P}_{\infty,\ldots,\infty}\bigl(T_h=T_h^i
\mbox{ for some }i\in\{1,\ldots ,N\}\bigr)=g(h),
\end{eqnarray*}
which completes the proof of (\ref{LBsame}).
\end{pf}

As a result of Proposition \ref{thm:false_alarm}, when $\mu_i=\mu$ for
all $i=1,\ldots,N$,
for any $\gamma>0$ and any
$N$-dimensional, predictable, nonsingular, stochastic instantaneous
correlation matrix $\Sigma_t$, we can choose the threshold $h$ using
\begin{equation}
\label{equalH}
\mathbb{E}_\infty^1\bigl\{T_h^1
\bigr\} \equiv \frac{2}{\mu^2}\bigl(\mathrm{e}^h-h-1\bigr)={N}
\gamma.
\end{equation}
Then we will have $T_\hbar\in\mathcal{T}_\gamma$. Moreover,
Proposition \ref{UpB} implies that, both the optimal detection delay
$\inf_{T\in\mathcal{T}_\gamma}J^{(N)}(T)$ and the detection delay of
this $N$-CUSUM
stopping rule $J^{(N)}(T_h)$, are bounded above by $\frac{2}{\mu^2} g(-h)$.

\subsubsection{Unequal drift case---complete information about \texorpdfstring{$\mu_i$}{mui}'s}
In this case, it is assumed that all $\mu_i$ are known and
$0 < \mu_1 =\mu_2 = \cdots= \mu_k < \min_{i>k} \mu_i$ holds for some
$k\in\{1,\ldots,N-1\}$. Then the monotonicity of function $g$ and
(\ref{hi}) imply that
$h_1 = h_2 = \cdots= h_k$, for the $N$-CUSUM stopping rule $T_{\hbar}$.
When $h_i$'s are all big, the condition (\ref{hi}) is approximately a
linear constraint on $h_i$'s, and hence $h_1<\min_{i>k} h_i$.
Intuitively, with high chances, $T_{h_i}^i$ for $1\le i\le k$ will
proceed $T_{h_j}^j$ for $k+1\le j\le N$ due to their smaller
thresholds. Hence, it is expected that $\mathbb{E}_{\infty,\ldots
,\infty}\{
T_\hbar\}\approx\mathbb{E}_{\infty,\ldots,\infty}\{
T_{h_1}^1\wedge\cdots\wedge
T_{h_k}^k\}\ge\frac{2}{k\mu_1^2}g(h_1)$, where the inequality follows
from (\ref{LBsame}). Below we rigorously show the validness of this
heuristic argument.


%
\begin{prop}\label{thm:false_alarm_k}
Suppose that the drifts $\mu_i$ of the observation processes $\xi^{(i)}_t$,
$i=1,\ldots,N$ are such that $0 < \mu_1 = \mu_2 = \cdots= \mu_k <
\min_{i>k} \mu_i$ holds.
Suppose also that the thresholds $\hbar$ satisfy (\ref{hi}).
Then\vadjust{\goodbreak} the mean time to the first false alarm for the $N$-CUSUM stopping rule
$T_\hbar$ satisfies
%
\begin{eqnarray}
\mathbb{E}_{\infty,\ldots,\infty}\{ T_\hbar\} &\ge &   \Biggl(1-
\sum_{j=k+1}^N\frac
{\mathbb{E}_\infty^1\{T_{h_1}^1\}}{\mathbb{E}_\infty^j\{
T_{h_{j}}^j\}} \Biggr)
\frac
{1}{k} \mathbb{E}_\infty^1\bigl
\{T_{h_1}^1\bigr\}
\nonumber
\\[-8pt]
\label{LBsamek}
\\[-8pt]
\nonumber
&=&  \Biggl(1-\sum
_{j=k+1}^N\frac{\mu
_j^2}{\mu_1^2}\frac{g(h_1)}{g(h_{j})}
\Biggr) \frac{2}{k\mu_1^2} g(h_1),
\end{eqnarray}
where the function $g$ is defined in (\ref{g}).
\end{prop}

\begin{pf}
Let us denote by $R_{h_1}:=T_{h_1}^{1} \wedge\cdots\wedge T_{h_1}^{k}$.
For any $k+1\le j\le N$,
following similar arguments to the ones in (\ref{R1}) through (\ref{R3}),
we have in this case that
%
\begin{eqnarray}
\frac{2}{\mu_j^2}g(h_j)&=& \mathbb{E}_{\infty,\ldots,\infty}
\bigl\{ T_{h_j}^j \bigr\}
\nonumber\\
&=&  \mathbb{E}_{\infty,\ldots,\infty}\{ T_\hbar\} + \mathbb{E}_{\infty,\ldots,\infty}
\bigl\{ \bigl(T_{h_j}^j - T_\hbar\bigr)
\mathbf{1}_{\{T_\hbar \neq T_{h_j}^j\}} \bigr\}
\nonumber
\\
\label{R12}
&= &\mathbb{E}_{\infty,\ldots,\infty}\{ T_\hbar\} + \mathbb{E}_{\infty,\ldots,\infty}
\bigl\{ \mathbb{E}_{\infty
,\ldots,\infty}\bigl\{ T_{h_j}^j -
T_\hbar | \mathcal{F}_{T_\hbar} \bigr\} \mathbf{1}_{\{T_\hbar \neq
T_{h_j}^j\}}
\bigr\}
\\
&\le& \mathbb{E}_{\infty,\ldots,\infty}\bigl\{T_{h_1}^1\bigr\}+
\frac{2}{\mu
_j^2} \mathbb{E} _{\infty,\ldots,\infty} \bigl\{ \bigl(g
\bigl(y_{T_{h_j}^j}^{(j)} \bigr)-g \bigl(y_{T_\hbar}^{(j)}
\bigr) \bigr) \mathbf{1}_{\{T_\hbar\neq T_{h_j}^j\}
} \bigr\}
\nonumber
\\
&\le& \frac{2}{\mu_1^2} g(h_1)+\frac{2}{\mu_j^2} g(h_j)
\mathbb{P} _{\infty,\ldots,\infty}\bigl(T_\hbar\neq T_{h_j}^j
\bigr),
\nonumber
\end{eqnarray}
which implies that
%
\begin{equation}
\label{eq:bound_prob}
\mathbb{P}_{\infty,\ldots,\infty}\bigl(T_\hbar=T_{h_j}^j
\bigr) \le\frac
{\mu_j^2}{\mu
_1^2} \frac{g(h_1)}{g(h_j)}.
\end{equation}
On the other hand, for any $1\le i\le k$, we similarly have
\[
\frac{2}{\mu_1^2} g(h_1)=\mathbb{E}_{\infty,\ldots,\infty}\bigl\{
T_{h_1}^{1}\bigr\} \le  \mathbb{E}_{\infty,\ldots,\infty}
\{T_\hbar\} + \frac{2}{\mu_1^2} g(h_1)
\mathbb{P}_{\infty,\ldots,\infty
}\bigl(T_\hbar\neq T_{h_1}^i
\bigr),
\]
which implies that
\[
\mathbb{E}_{\infty,\ldots,\infty}\{T_\hbar\} \ge\frac{2}{\mu
_1^2}
g(h_1) \mathbb{P}_{\infty,\ldots,\infty}\bigl(T_\hbar=T_{h_1}^i
\bigr).
\]
Summing up both sides of the above inequality for all $1\le i\le k$, we
obtain that
%
\begin{eqnarray}
k \mathbb{E}_{\infty,\ldots,\infty}\{T_\hbar\} &\ge & \frac{2}{\mu
_1^2}
g(h_1) \mathbb{P}_{\infty,\ldots,\infty}(T_\hbar=R_{h_1})
\nonumber
\\[-8pt]
\label{FAL2P}
\\[-8pt]
\nonumber
&=&
\frac{2}{\mu_1^2} g(h_1) \bigl[1-\mathbb{P}_{\infty,\ldots,\infty}(T_\hbar
\neq R_{h_1})\bigr].
\end{eqnarray}
However, we also have
%
\begin{equation}\label{P2g}
\hspace*{15pt}\mathbb{P}_{\infty,\ldots,\infty}(T_\hbar\neq R_{h_1})\le\sum
_{j=k+1}^N\mathbb{P} _{\infty,\ldots,\infty}
\bigl(T_\hbar=T_{h_j}^j\bigr)\le\sum
_{j=k+1}^N\frac
{\mu
_j^2}{\mu_1^2} \frac{g(h_1)}{g(h_j)},
\end{equation}
where we used (\ref{eq:bound_prob}) in the above inequality. It follows
from (\ref{FAL2P}) and (\ref{P2g}) that
%
%
\begin{equation}
\label{eq:bound_fa}
\mathbb{E}_{\infty,\ldots,\infty}\{T_\hbar\} \ge \Biggl(1-\sum
_{j=k+1}^N\frac{\mu_j^2}{\mu_1^2}
\frac
{g(h_1)}{g(h_j)} \Biggr) \frac{2}{k\mu_1^2}g(h_1),
\end{equation}
%
which completes the proof.
\end{pf}

As a result of Proposition \ref{thm:false_alarm_k}, when $\mu
_1=\cdots
=\mu_k <
\min_{i>k} \mu_i$, then for any $\gamma>0$ and any $N$-dimensional,
predictable, nonsingular, stochastic instantaneous correlation matrix
$\Sigma_t$,
we can choose the set of thresholds $\hbar$ using (\ref{hi}) and the
transcendental equation
%
\begin{eqnarray}
&& \Biggl(1-\sum_{j=k+1}^N
\frac{\mathbb{E}_\infty^1\{T_{h_1}^1\}
}{\mathbb{E}_\infty^j\{
T_{h_{j}}^j\}} \Biggr) \frac{1}{k} \mathbb{E}_\infty^1
\bigl\{T_{h_1}^1\bigr\}
\nonumber
\\[-8pt]
\label{inequalH}
\\[-8pt]
\nonumber
&& \qquad = \Biggl(1-\sum
_{j=k+1}^N\frac{\mu_j^2}{\mu_1^2}\frac{\mathrm{e}^{h_1}-h_1-1}{\mathrm{e}^{h_{j}}-h_{j}-1}
\Biggr)\frac{2}{k\mu
_1^2}\bigl(\mathrm{e}^{h_1}-h_1-1\bigr)=
\gamma,
\end{eqnarray}
%
then the resulting $N$-CUSUM stopping rule $T_\hbar\in\mathcal
{T}_\gamma$.
Again, Proposition \ref{UpB} then implies that, both the optimal
detection delay $\inf_{T\in\mathcal{T}_\gamma}J^{(N)}(T)$ and the
detection delay of this $N$-CUSUM
stopping rule $J^{(N)}(T_\hbar)$, are bounded above by
$\frac{2}{\mu_1^2}g(-h_1)$.

\subsubsection{The general case---partial information about \texorpdfstring{$\mu_i$}{mui}'s}
In this case, it is assumed that only $
\mu_1, \underline{\mu}_i, \overline{\mu}_i$, $i=2,\ldots,N$, are known,
and that $0<\mu_1\le\underline{\mu}_i\le\mu_i\le\overline{\mu}_i$.
Without loss of generality, we assume that $0<\mu_1=\underline{\mu
}_2=\cdots=\underline{\mu}_{k'}<\min_{i>k'}\underline{\mu}_i$
holds for
some $k'=\{1,\ldots,N-1\}$.
%

\begin{prop}\label{thm:false_alarm_k1}
Suppose that the drifts $\mu_i$ of the observation processes $\xi^{(i)}_t$,
$i=1,\ldots,N$ are such that $0 < \mu_1 = \underline{\mu}_2= \cdots=
\underline{\mu}_{k'} <
\min_{i>k'} \underline{\mu_i}$ holds\vspace*{1pt} and $\mu_i\in[\underline{\mu
}_i,\overline{\mu}_i]$ for all $i=2,\ldots,N$.\vspace*{1pt}
Suppose also that the thresholds $\hbar$ satisfy (\ref{hi1}).
Then the mean time to the first false alarm for the $N$-CUSUM stopping rule
$T_\hbar$ satisfies
%
\begin{eqnarray}
&& \mathbb{E}_{\infty,\ldots,\infty}\{ T_\hbar\}
\nonumber
\\[-4pt]
\label{LBsamek1}\\[-12pt]
\nonumber
&& \qquad \ge
\frac{2}{\sum_{1\le i\le
k'}\mu_1(2\overline{\mu}_i-\mu_1)} \biggl(1-\sum_{k'+1\le j\le
N}
\frac
{\underline{\mu}_j(2\overline{\mu}_j-\underline{\mu_j})}{\mu_1^2} \frac{g(h_1)}{g(h_j)} \biggr) g(h_1),\hspace*{-23pt}
\end{eqnarray}
where the function $g$ is defined in (\ref{g}).
\end{prop}
\begin{pf}
According to (\ref{hi1}), we have $h_1=h_2=\cdots=h_{k'}$. Let us
denote by $R_{h_1}:=T_{h_1}^{1} \wedge\cdots\wedge T_{h_1}^{k'}$.
Similar (\ref{eq:bound_prob}) in the proof of Proposition \ref
{thm:false_alarm_k1}, for any $k'+1\le j\le N$, we have
\begin{eqnarray*}
g(h_j) &=& \mathbb{E}_{\infty,\ldots,\infty} \biggl\{\int
_0^{T_{h_j}^j}\underline {\mu}_j \biggl(
\mu_j-\frac{1}{2}\underline{\mu}_j \biggr)\, d s \biggr\}
\\
&=&  \mathbb{E}_{\infty,\ldots,\infty} \biggl\{\int_0^{T_\hbar
}
\underline{\mu }_j \biggl(\mu_j-\frac{1}{2}
\underline{\mu}_j \biggr)\,\dd s \biggr\}
\\
&&{}+ \mathbb{E}_{\infty,\ldots,\infty} \biggl\{\mathbb{E}_{\infty
,\ldots,\infty} \biggl\{ \int
_{T_\hbar}^{T_{h_j}^j}\underline{\mu}_j \biggl(
\mu_j-\frac
{1}{2}\underline{\mu}_j \biggr)\,\dd s \Big|
\mathcal{F}_{T_\hbar} \biggr\}\mathbf{1}_{\{T_\hbar\neq
T_{h_j}^j\}} \biggr\}
\\
&\le & \mathbb{E}_{\infty,\ldots,\infty}
\biggl\{\int_0^{T_{h_1}^1}\underline{\mu
}_j \biggl(\mu_j-\frac{1}{2}\underline{
\mu}_j \biggr)\,\dd s \biggr\}\\
&&{}+ \mathbb{E}_{\infty,\ldots,\infty} \bigl\{
\bigl(g \bigl(y_{T_{h_j}^j}^{(j)} \bigr)-g \bigl(y_{T_\hbar}^{(j)}
\bigr) \bigr) \mathbf{1}_{\{T_\hbar\neq
T_{h_j}^j\}} \bigr\}
\\
&\le & \frac{\underline{\mu}_j(2\overline{\mu}_j-\underline{\mu
_j})}{2}\mathbb{E} _{\infty,\ldots,\infty}\bigl\{T_{h_1}^1
\bigr\}+ g(h_j) \mathbb{P}_{\infty
,\ldots
,\infty}\bigl(T_\hbar\neq
T_{h_j}^j\bigr)
\\
&=& \frac{\underline{\mu}_j(2\overline{\mu}_j-\underline{\mu
_j})}{\mu
_1^2}g(h_1)+ g(h_j) \mathbb{P}_{\infty,\ldots,\infty}
\bigl(T_\hbar\neq T_{h_j}^j\bigr).
\end{eqnarray*}
It follows that
\[
\mathbb{P}_{\infty,\ldots,\infty}\bigl(T_\hbar=T_{h_j}^j
\bigr) \le\frac
{\underline{\mu
}_j(2\overline{\mu}_j-\underline{\mu_j})}{\mu_1^2} \frac
{g(h_1)}{g(h_j)},
\]
which implies that
%
\begin{eqnarray}
\qquad\quad\mathbb{P}_{\infty,\ldots,\infty}(T_\hbar\neq R_{h_1})
&=& \mathbb{P}_{\infty,\ldots,\infty}\bigl(T_\hbar=T_{h_j}^j,
\mbox{for some } j \in\bigl\{k'+1, \ldots, N\bigr\} \bigr)
\nonumber
\\[-8pt]
\label{eq:bound_prob1}
\\[-8pt]
\nonumber
&\le & \sum_{k'+1\le j\le N}\frac{\underline{\mu}_j(2\overline{\mu
}_j-\underline{\mu_j})}{\mu_1^2}
\frac{g(h_1)}{g(h_j)}.
\end{eqnarray}
On the other hand, for any $1\le i\le k'$, by It\^{o}'s formula
(see, e.g., \cite{Okse}, Theorem~4.1.2) we have
\begin{eqnarray*}
g(h_1) &=& \mathbb{E}_{\infty,\ldots,\infty} \biggl\{\int_0^{T_{h_1}^i}
\mu_1 \biggl(\mu_i-\frac{1}{2}\mu_1
\biggr) \,\dd s \biggr\}
\nonumber
\\
&=&  \mathbb{E}_{\infty,\ldots,\infty} \biggl\{ \int_0^{T_\hbar}
\mu _1 \biggl(\mu _i-\frac{1}{2}\mu_1
\biggr) \,\dd s \biggr\}
\nonumber
\\
&&{}+ \mathbb{E}_{\infty,\ldots,\infty} \biggl\{\mathbb{E}_{\infty
,\ldots,\infty} \biggl\{ \int
_{T_\hbar}^{T_{h_1}^i} \mu_1 \biggl(
\mu_i-\frac{1}{2}\mu _1 \biggr) \,\dd s \Big|
\mathcal{F}_{T_\hbar} \biggr\} \mathbf{1}_{\{T_\hbar\neq
T_{h_1}^i\}} \biggr\}
\nonumber
\\
&=& \mathbb{E}_{\infty,\ldots,\infty} \biggl\{ \int_0^{T_\hbar}
\mu _1 \biggl(\mu _i-\frac{1}{2}\mu_1
\biggr) \,\dd s \biggr\}\\
&&{}+ \mathbb{E}_{\infty,\ldots,\infty} \bigl\{ \bigl(g
\bigl(y_{T_{h_1}^i}^{(i)} \bigr)-g \bigl(y_{T_\hbar}^{(i)}
\bigr) \bigr)\mathbf{1}_{\{T_\hbar\neq
T_{h_1}^i\}} \bigr\}
\nonumber
\\
&\le & \frac{\mu_1(2\overline{\mu}_i-\mu_1)}{2}\mathbb{E}_{\infty
,\ldots,\infty
}\{T_\hbar\}+
g(h_1) \mathbb{P}_{\infty,\ldots,\infty}\bigl(T_\hbar \neq
T_{h_1}^i\bigr),
\nonumber
\end{eqnarray*}
which implies that
\[
\frac{\mu_1(2\overline{\mu}_i-\mu_1)}{2}\mathbb{E}_{\infty
,\ldots,\infty}\{ T_\hbar\}\ge
g(h_1)\mathbb{P}_{\infty,\ldots,\infty}\bigl(T_\hbar
=T_{h_1}^i\bigr).
\]
Summing up both sides of the above inequality for all $1\le i\le k'$,
we obtain that
\begin{eqnarray}
&& \frac{\sum_{1\le i\le k'}\mu_1(2\overline{\mu}_i-\mu
_1)}{2}\mathbb{E}_{\infty
,\ldots,\infty}\{T_\hbar\}\nonumber \\
&& \qquad \ge
g(h_1)\mathbb{P}_{\infty,\ldots
,\infty}(T_\hbar =R_{h_1})
\nonumber
\\[-8pt]
\\[-8pt]
\nonumber
&&\qquad=g(h_1)\bigl[1-\mathbb{P}_{\infty,\ldots,\infty}(T_\hbar\neq
R_{h_1})\bigr]
\\
&&\qquad\ge g(h_1) \biggl(1-\sum_{k'+1\le j\le N}
\frac{\underline{\mu
}_j(2\overline{\mu}_j-\underline{\mu_j})}{\mu_1^2} \frac
{g(h_1)}{g(h_j)} \biggr),
\nonumber
\end{eqnarray}
where we used (\ref{eq:bound_prob1}) in the last step. The conclusion
of the proposition follows immediately.
\end{pf}
As a result of Proposition \ref{thm:false_alarm_k1}, when we only known
$\mu_1$ and possible ranges for other drift $\mu_i$'s, given any
$\gamma
>0$ and any $N$-dimensional, predictable, nonsingular, stochastic
instantaneous correlation matrix $\Sigma_t$,
we can choose the set of thresholds~$\hbar$ using (\ref{hi1}) and the
transcendental equation
%
\begin{equation}
\label{inequalH1}
\qquad\hspace*{5pt}\Biggl(1-\sum_{j=k'+1}^N
\frac{\underline{\mu}_j(2\overline{\mu
}_j-\underline{\mu}_j)}{\mu_1^2}\frac{\mathrm{e}^{h_1}-h_1-1}{\mathrm{e}^{h_{j}}-h_{j}-1} \Biggr)\frac{2(\mathrm{e}^{h_1}-h_1-1)}{\sum_{1\le
i\le k'}\mu_1(2\overline{\mu}_i-\mu_1)}=\gamma,
\end{equation}
%
then the resulting $N$-CUSUM stopping rule $T_\hbar\in\mathcal
{T}_\gamma
$. Again, Proposition \ref{UpB} then implies that both the optimal
detection delay $\inf_{T\in\mathcal{T}_\gamma}J^{(N)}(T)$ and the
detection delay of this $N$-CUSUM
stopping rule $J^{(N)}(T_\hbar)$, are bounded above by
$\frac{2}{\mu_1^2}g(-h_1)$.

\subsection{The lower bound}
In this subsection, we present a robust lower bound for the optimal
detection delay $\inf_{T\in\mathcal{T}_\gamma}J^{(N)}(T)$. In fact, we
can prove a stronger statement: for any stopping rule $T\in\mathcal
{T}_\gamma$, its detection delay $J^{(N)}(T)$, is bounded below by the
optimal detection delay in one dimension. The proof is accomplished by
a change of measure argument as in \cite{Mous04} plus a decomposition
formula for the Radon--Nikodym derivative in $N$ dimensions.

%
\begin{lem}\label{lem:Girsanov}
Let $\mathbb{Q}=\mathbb{P}_{\infty,\ldots,\infty}$ be the law of
the $N$-tuple
process $W_t:=(w_t^{(1)},w_t^{(2)},\ldots,w_t^{(N)})$ for the Brownian
motions defined in (\ref{Itodynamics}). And let $\mathbb{Q}_1$ be the
law of
the $N$-tuple process $(\mu_1t+w_t^{(1)}, w_t^{(2)},\ldots,
w_t^{(N)})$. Then for all $t>0$,
%
\begin{equation}
\frac{d\mathbb{Q}_1}{d\mathbb{Q}} \bigg|_{\mathcal{F}_t}=\mathrm{e}^{u_t^{(1)}}\cdot\mathcal {E}
\bigl(B^{(1)}\bigr)_t,
\end{equation}
where $u_\cdot^{(1)}$ is defined in (\ref{1dimpr}) and
$\mathcal{E}(B^{(1)})_{\cdot}$
is the stochastic exponential of the local martingale $B^{(1)}_\cdot$
defined in
(\ref{Bt})--(\ref{wtB}).
Moreover, the standard Brownian motions driving $B^{(1)}_\cdot$ are
independent of $w^{(1)}$.
\end{lem}

\begin{pf}
The proof can be found in the \hyperref[pf]{Appendix}.
\end{pf}
%

\begin{prop} \label{lemmaLB}
For any stopping rule $T\in\mathcal{T}_\gamma$, we have
$J^{(N)}(T) \ge(2 / \mu_1^2) g(-\nu^\star_1)$,
where\vspace*{1pt} $\nu^\star_i$ satisfies $g(\nu^\star_1)=(\mu_1^2 / 2) \gamma$
for the
function $g$ defined in~(\ref{g}).
\end{prop}

\begin{pf}
%
Let $T$ be an arbitrary $\mathbb{F}$-stopping rule such that
$\mathbb{E}_{\infty,\ldots,\infty} \{ T \} \ge\gamma$
holds and observe that
%
\begin{equation}
\label{new index}
\qquad \quad J^{(N)}(T) \geq\tilde{J}^{(N)}_1(T)
:= \mathop{\sup_{s\in\mathbb
{R}_+}} \operatorname{essup} \mathbb{E}_{s,\infty,\ldots,\infty}
\bigl\{ (T - s)^+ |\mathcal {F}_{s} \bigr\} \geq
\tilde{J}^{(N)}_1(T_\nu),
\end{equation}
where $T_\nu:= T \wedge T_\nu^{1} \le T$, a.s., and $T_\nu^{1}$ is
the CUSUM
stopping rule given in (\ref{1dimCUSUMstop}) for some threshold $\nu$
which will be determined later. Clearly, $T_\nu$ is a finite stopping rule.
In what follows, we will demonstrate that for any given $\varepsilon >0$, there
exists a $\nu>0$ such that
%
\begin{equation}
\label{INEQ2}
\tilde{J}^{(N)}_1(T_\nu) \ge
\frac{2}{\mu_1^2} g\bigl(-\nu^*_1\bigr) - \varepsilon,
\end{equation}
where $\nu_1^\star$ is chosen so that $g(\nu_1^\star)=(\mu
_1^2/2)\gamma$ and the function $g$ is given by (\ref{g}). Because $\varepsilon$ in
(\ref{INEQ2}) can be arbitrarily small, (\ref{new index}) and (\ref
{INEQ2}) will imply the assertion in the proposition for $i=1$. This is
in a similar light as in \cite{Mous04}.

By\vspace*{1pt} applying It\^{o}'s formula
(see, e.g., \cite{Okse}, Theorem 4.1.2 and \cite{Mous04})
to  $\{g(-y_t^{(1)})\}_{s \wedge T_\nu \le
t < T_\nu}$ and proceed by using similar arguments as in
(\ref{Ito's expansion})--(\ref{Ito}) in Proposition \ref{thm:false_alarm},
we obtain that, for any fixed $s\in\mathbb{R}_+$,
%
\begin{eqnarray}
&& \mathbb{E}_{s,\infty,\ldots,\infty} \bigl\{ (T_{\nu} - s)^+ |
\mathcal{F}_{s} \bigr\}
\nonumber
\\[-8pt]
\label{Ito 1d}
\\[-8pt]
\nonumber
&& \qquad  = \frac{2}{\mu_1^2} \mathbb{E}_{s,\infty,\ldots,\infty}
\bigl\{ g\bigl(-y_{T_\nu
}^{(1)}\bigr) - g\bigl(-y_{s}^{(1)}
\bigr) | \mathcal{F}_{s}\bigr\} {\bf1}_{\{T_\nu
\ge s \}}.\\[-18pt]
\nonumber
\end{eqnarray}\eject
\noindent Using Girsanov's theorem (see, e.g., \cite{KS}, Chapter~3, Theorem~5.1)
and Lemma~\ref{lem:Girsanov} at the finite stopping rule $T_\nu\wedge
n$ for a fixed $n>0$, we have that
\begin{eqnarray*}
&& \mathbb{E}_{s,\infty,\ldots,\infty}\bigl\{g\bigl(-y_{T_\nu\wedge n}^{(1)}\bigr)
- g\bigl(-y_{s}^{(1)}\bigr) | \mathcal{F}_{s}\bigr
\} \mathbf{1}_{\{T_\nu\wedge n>s\}
}
\\
&& \qquad = \mathbb{E}_{\infty,\ldots,\infty} \biggl\{ \mathrm{e}^{u_{T_\nu
\wedge n}^{(1)}
- u_{s}^{(1)}} \cdot
\frac{\mathcal{E}(B^{(1)})_{T_\nu\wedge
n}}{\mathcal{E}(B^{(1)})_{s}} \bigl[g\bigl(-y_{T_\nu\wedge
n}^{(1)}\bigr)-g
\bigl(-y_{s}^{(1)}\bigr)\bigr] \Big| \mathcal{F}_{s}
\biggr\}\\
&& \quad\qquad{}\times\mathbf{1}_{\{T_\nu
\wedge n>s\}}.
\end{eqnarray*}
Consider the enlargement of filtration $\mathcal{F}_t^a=\mathcal
{F}_t\vee\mathcal{G}_{T_\nu^1}^{(1)}$. Then clearly, $\mathcal
{F}_t^a\equiv\mathcal{F}_t$ for all $t\ge T_\nu^1$, but on the event
$\{
T_\nu\wedge n>s\}$, for all $t\in[s, T_\nu\wedge n]\subset[0,T_\nu^1]$,
we have $\mathcal{F}_t\subsetneq\mathcal{F}_t^a$. By the tower property
of conditional expectation, on the event that $\{T_\nu\wedge n>s\}$,
%
\begin{eqnarray}
&& \mathbb{E}_{\infty,\ldots,\infty} \biggl\{ \mathrm{e}^{u_{T_\nu
\wedge n}^{(1)}
- u_{s}^{(1)}}
\cdot\frac{\mathcal{E}(B^{(1)})_{T_\nu\wedge
n}}{\mathcal{E}(B^{(1)})_{s}} \bigl[g\bigl(-y_{T_\nu\wedge
n}^{(1)}\bigr)-g
\bigl(-y_{s}^{(1)}\bigr)\bigr] \Big| \mathcal{F}_{s}
\biggr\}
\nonumber
\\
&&\qquad =\mathbb{E}_{\infty,\ldots,\infty} \biggl\{\mathrm{e}^{u_{T_\nu
\wedge n}^{(1)}
- u_{s}^{(1)}} \bigl[g
\bigl(-y_{T_\nu\wedge n}^{(1)}\bigr)-g\bigl(-y_{s}^{(1)}
\bigr)\bigr]
\nonumber
\\[-8pt]
\label{Girsanov1}\\[-8pt]
\nonumber
&&{}\hspace*{47pt}\qquad\quad{}\times \mathbb{E} _{\infty,\ldots,\infty} \biggl\{\frac{\mathcal{E}(B^{(1)})_{T_\nu
\wedge
n}}{\mathcal{E}(B^{(1)})_{s}} \Big|
\mathcal{F}_{s}^a \biggr\} \Big| \mathcal {F}_{s}
\biggr\}
\\
&& \qquad =\mathbb{E}_{\infty,\ldots,\infty} \bigl\{ \mathrm{e}^{u_{T_\nu
\wedge n}^{(1)}
- u_{s}^{(1)}} \bigl[g
\bigl(-y_{T_\nu\wedge n}^{(1)}\bigr)-g\bigl(-y_{s}^{(1)}
\bigr)\bigr] | \mathcal{F}_{s} \bigr\},
\nonumber
\end{eqnarray}
where the last equality is due to the fact that $\mathcal
{E}(B^{(1)})_{t}$ is a $\overline{\mathbb{F}}$-exponential martingale [under
assumption equation (\ref{eq:Nov})] driven by Brownian motions that are
independent of $w^{(1)}$ (see Lemma
\ref{lem:Girsanov}). Similarly, it can be shown that
%
\begin{equation}
\label{Girsanov2}
\mathbf{1}_{\{T_\nu\wedge n>s\}} = \mathbf{1}_{\{T_\nu\wedge n>s\}}
\mathbb{E}_{\infty,\ldots,\infty}\bigl\{\mathrm{e}^{u_{T_\nu\wedge n}^{(1)}-u_{s}^{(1)}} |
\mathcal{F}_{s}\bigr\}.
\end{equation}
We now let $n\uparrow\infty$ in (\ref{Girsanov1}) and (\ref{Girsanov2}). From the fact that 
$u_{T_\nu\wedge n}^{(1)} - u_{s}^{(1)} \le u_{T_\nu\wedge n}^{(1)} -
m_{T_\nu\wedge n}^{(1)} = y_{T_\nu\wedge n}^{(1)} \le\nu$, and the
monotonicity of function $g(-h)$,
we have that
\begin{eqnarray}
0&<& \mathrm{e}^{u_{T_\nu\wedge n}^{(1)} - u_{s}^{(1)}} \le\mathrm{e}^\nu < \infty \quad  \mbox{and}\quad
\bigl|g\bigl(-y_{T_\nu\wedge n}^{(1)}\bigr) - g\bigl(-y_{s}^{(1)}
\bigr) \bigr| \le2 g(-\nu) <\infty\nonumber\\
\eqntext{\forall n>0}
\end{eqnarray}
and thus the bounded convergence theorem implies that, on the event
$\{T_\nu>s\}$,
%
\begin{eqnarray}
\label{E88}
&& 1 =  \mathbb{E}_{\infty,\ldots,\infty}\bigl\{\mathrm{e}^{u_{T_\nu
}^{(1)}-u_{s}^{(1)}} |
\mathcal{F}_{s}\bigr\},
\\
&& \mathbb{E}_{s,\infty,\ldots,\infty}\bigl\{g\bigl(-y_{T_\nu}^{(1)}
\bigr)-g\bigl(-y_{s}^{(1)}\bigr) | \mathcal{F}_{s}
\bigr\}
\nonumber
\\[-8pt]
\label{E88g}\\[-8pt]
\nonumber
&&\qquad  = \mathbb{E}_{\infty,\ldots,\infty}\bigl\{\mathrm{e}^{u_{T_\nu
}^{(1)}-u_{s}^{(1)}} \bigl[g
\bigl(-y_{T_\nu}^{(1)}\bigr)-g\bigl(-y_{s}^{(1)}
\bigr)\bigr] | \mathcal{F}_{s}\bigr\}.
\end{eqnarray}
It follows from (\ref{new index}), (\ref{Ito 1d}), (\ref{E88}) and
(\ref{E88g}) that
%
\begin{eqnarray}
&& \tilde{J}_1^{(N)}(T_\nu)
\mathbb{E}_{\infty,\ldots,\infty}\bigl\{ \mathrm{e}^{u_{T_\nu}^{(1)}-u_{s}^{(1)}} |
\mathcal{F}_{s}\bigr\}\mathbf{1}_{\{
T_\nu >s\}}\nonumber \\
&&\qquad = \tilde{J}_1^{(N)}(T_\nu)
\mathbf{1}_{\{T_\nu>s\}}\nonumber
\\[-8pt]
\label{ineq}
\\[-8pt]
\nonumber
&& \qquad \ge \mathbb{E}_{s,\infty,\ldots,\infty}\bigl\{(T_\nu-s)^+ | \mathcal{F}_{s}\bigr\}\mathbf{1}_{\{T_\nu>s\}}
\\
&&\qquad= \frac{2}{\mu_1^2}\mathbb{E}_{\infty,\ldots,\infty}\bigl\{\mathrm{e}^{u_{T_\nu
}^{(1)}-u_{s}^{(1)}}
\bigl[g\bigl(-y_{T_\nu}^{(1)}\bigr)-g\bigl(-y_{s}^{(1)}
\bigr)\bigr] | \mathcal{F}_{s}\bigr\}\mathbf{1}_{\{T_\nu>s\}}.
\nonumber
\end{eqnarray}
Following the same arguments as in Theorem 2 of \cite{Mous04}, we
integrate both sides of the above inequality with respect to $(-\dd
m_{s}^{(1)})$ for all $s \in[0,T_\nu]$ and then take the expectation
under $\mathbb{P}_{\infty,\ldots,\infty}$.
We therefore obtain that
\begin{eqnarray*}
&& \tilde{J}_1^{(N)}(T_\nu) \mathbb{E}_{\infty,\ldots,\infty}
\biggl\{ \mathrm{e}^{u_{T_\nu}^{(1)}}\int_0^{T_\nu}
\mathrm{e}^{-u_{s}^{(1)}}\bigl(-\dd m_{s}^{(1)}\bigr) \biggr\}
\\
&& \qquad \ge\frac{2}{\mu_1^2} \mathbb{E}_{\infty,\ldots,\infty} \biggl\{ \mathrm{e}^{u_{T_{\nu}}^{(1)}}\int_0^{T_\nu}\mathrm{e}^{-u_{s}^{(1)}}\bigl[g\bigl(-y_{T_\nu}^{(1)}\bigr)-g
\bigl(-y_{s}^{(1)}\bigr)\bigr]\bigl(-\dd m_{s}^{(1)}
\bigr) \biggr\}.
\end{eqnarray*}
Notice that the measure $\dd m_{s}^{(1)}$ is supported on the random
set $\{s \mid y_{s}^{(1)}=0 \} = \{s \mid u_{s}^{(1)} = m_{s}^{(1)} \}
$, and that $g(0)=0$, thus we obtain that
\[
\tilde{J}_1^{(N)}(T_\nu) \mathbb{E}_{\infty,\ldots,\infty}
\bigl\{ \mathrm{e}^{y_{T_\nu}^{(1)}} - \mathrm{e}^{u_{T_\nu}^{(1)}}\bigr\} \ge
\frac{2}{\mu_1^2} \mathbb{E}_{\infty,\ldots,\infty}\bigl\{\bigl[\mathrm{e}^{y_{T_\nu}^{(1)}} - \mathrm{e}^{u_{T_\nu}^{(1)}}\bigr] g\bigl(-y_{T_\nu
}^{(1)}
\bigr)\bigr\}.
\]
On the other hand, by letting $s = 0$ in (\ref{ineq}) we have that
\[
\tilde{J}_1^{(N)}(T_\nu) \mathbb{E}_{\infty,\ldots,\infty}
\bigl\{ \mathrm{e}^{u_{T_\nu}^{(1)}}\bigr\} \ge\frac{2}{\mu_1^2}
\mathbb{E}_{\infty,\ldots,\infty}\bigl\{\mathrm{e}^{u_{T_\nu}^{(1)}} g\bigl(-y_{T_\nu}^{(1)}
\bigr)\bigr\}.
\]
In all, we have that
%
\begin{equation}
\label{INEQ3}
\tilde{J}_1^{(N)}(T_\nu) {
\mathbb{E}_{\infty,\ldots,\infty}\bigl\{ \mathrm{e}^{y_{T_\nu}^{(1)}}\bigr\}} \ge
\frac{2}{\mu_1^2} {\mathbb{E}_{\infty,\ldots,\infty}\bigl\{\mathrm{e}^{y_{T_\nu}^{(1)}} g
\bigl(-y_{T_\nu
}^{(1)}\bigr)\bigr\}}
\end{equation}
holds.

To relate the detection delay in (\ref{INEQ3}) to the first false alarm
constraint $\gamma$, we use similar arguments as in (\ref{Ito's
expansion})--(\ref{Ito}) in
Proposition \ref{thm:false_alarm}, to obtain that
%
\begin{equation}
\label{flase rate}
\frac{2}{\mu_1^2} \mathbb{E}_{\infty,\ldots,\infty}\bigl\{ g
\bigl(y_{T_\nu
}^{(1)}\bigr) \bigr\} = \mathbb{E}_{\infty,\ldots,\infty}
\{ T_\nu\}.
\end{equation}
By taking the limit as $\nu\uparrow\infty$ and using monotone convergence
theorem, we have that $T_\nu\uparrow T$, and $\lim_{\nu\uparrow
\infty}
\mathbb{E}_{\infty,\ldots,\infty}\{T_\nu\} = \mathbb{E}_{\infty
,\ldots,\infty}\{T\}
\geq\gamma$, which implies that there exists a large enough $\nu$, such that
\[
\frac{2}{\mu_1^2} \mathbb{E}_{\infty,\ldots,\infty}\bigl\{g\bigl(y_{T_\nu
}^{(1)}
\bigr)\bigr\} = \mathbb{E}_{\infty,\ldots,\infty}\{T_\nu\} \ge \gamma-\varepsilon
\]
holds for any prespecified $\varepsilon >0$.
Now consider the nonnegative function $p(y) := \mathrm{e}^{y}
[g(-y) - g(-\nu^\star_1)] - g(y) + g(\nu^\star_1)$,
using which we trivially have $\mathbb{E}_{\infty,\ldots,\infty}\{
p(y_{T_\nu
}^{(1)}) \} \ge0$, implying that
%
\begin{eqnarray}
\mathbb{E}_{\infty,\ldots,\infty}\bigl\{\mathrm{e}^{y_{T_\nu}^{(1)}} g
\bigl(-y_{T_\nu
}^{(1)}\bigr)\bigr\} &\ge & \mathbb{E}_{\infty,\ldots,\infty}
\bigl\{\mathrm{e}^{y_{T_\nu
}^{(1)}}\bigr\} g\bigl(-\nu^\star_1
\bigr) + \mathbb{E}_{\infty,\ldots,\infty}\bigl\{g\bigl(y_{T_\nu}^{(1)}
\bigr)\bigr\} -g\bigl(\nu ^\star_1\bigr)
\nonumber\\
&=&  \mathbb{E}_{\infty,\ldots,\infty}\bigl\{\mathrm{e}^{y_{T_\nu
}^{(1)}}\bigr\} g\bigl(-\nu
^\star_1\bigr) + \mathbb{E}_{\infty,\ldots,\infty}\bigl\{g
\bigl(y_{T_\nu}^{(1)}\bigr)\bigr\} - \frac
{\mu_1^2}{2} \gamma
\nonumber
\\[-8pt]
\label{ineqe}
\\[-8pt]
\nonumber
&\ge & \mathbb{E}_{\infty,\ldots,\infty}\bigl\{\mathrm{e}^{y_{T_\nu
}^{(1)}}\bigr\} g\bigl(-
\nu^\star_1\bigr) - \frac{\mu_1^2}{2} \varepsilon
\\
&\ge & \mathbb{E}_{\infty,\ldots,\infty}\bigl\{\mathrm{e}^{y_{T_\nu
}^{(1)}}\bigr\} \biggl[g
\bigl(-\nu^\star_1\bigr) - \frac{\mu_1^2}{2} \varepsilon
\biggr],
\nonumber
\end{eqnarray}
since $\mathbb{E}_{\infty,\ldots,\infty}\{\mathrm{e}^{y^{(1)}_{T_\nu}}\}\ge1$.
The above inequality in (\ref{ineqe}) together with (\ref{INEQ3})
yields~(\ref{INEQ2}), which completes the proof.
\end{pf}

\section{Asymptotic optimality of the $N$-CUSUM stopping rule} \label{Asympt}

In this section, we demonstrate the asymptotic optimality of the
$N$-CUSUM stopping rule $T_{\hbar}$ for $\hbar$ chosen such that
(\ref{hi}) and either (\ref{equalH}) or (\ref{inequalH}) hold, or
(\ref{hi1}) and (\ref{inequalH1}) hold.
To this end, we examine the asymptotic behavior of the robust upper and
low bounds established in Section~\ref{Robust}. We show that the additional
detection delay of $T_\hbar$ over the optimal detection delay remains
bounded as the mean time of the first false alarm~$\gamma$ increases
without bound.

Let any sufficiently large $\gamma>0$ and recall from Section~\ref{Robust} (in particular Propositions \ref{UpB} and \ref{lemmaLB}) that
the optimal detection delay in (\ref{eqnproblem}) is bounded from below
and above as
%
\begin{equation}
\label{UBJ}
\frac{2}{\mu_1^2}g\bigl(-\nu_1^\star\bigr)
\le\inf_{T\in\mathcal
{T}_\gamma} J^{(N)}(T) \leq J^{(N)}(T_\hbar)
\le\frac{2}{\mu_1^2}g(-h_1),
\end{equation}
where the set of thresholds $\nu_1^\star$ and $\hbar$ is, respectively,
determined using $(\mu_1^2/2)\times\break g(\nu_1^\star)=\gamma$ and
either (\ref{hi}) together with (\ref{equalH}) or with (\ref{inequalH}), or (\ref{hi1}) together with (\ref{inequalH1}),
when the drifts' sizes $\mu_i$ 
are all known and equal or unequal,\footnote{Note that when the drifts
$\mu_i$ are different, we do not necessarily require the uniqueness of
$\hbar$ that solves (\ref{hi}) and (\ref{inequalH}).} or partially known,
respectively.

It is easily seen from Result 3 in the Appendix of \cite{HadjZhanPoor}
that, as $\gamma\to\infty$,
%
\begin{eqnarray}
\nu_1^\star&=&\log\frac{\mu_1^2}{2}+\log\gamma+o(1),
\\
\label{lowbound1}
\frac{2}{\mu_1^2}g\bigl(-\nu_1^\star\bigr) &=&
\frac{2}{\mu_1^2} \biggl(\log\frac{\mu_1^2}{2}+\log\gamma -1+o(1)
\biggr).
\end{eqnarray}

Moreover, when all the drifts are known and $\mu_i = \mu>0$ for all
$i=1,\ldots,N$, the thresholds $h_i=h>0$ for all $i=1,\ldots,N$. Using
(\ref{equalH}) and Result 3 in the Appendix of \cite{HadjZhanPoor} we
have that, as $\gamma\to\infty$,
%
\begin{eqnarray}
h_1&=&\log\frac{N \mu_1^2}{2}+\log\gamma+o(1),
\\[-2pt]
\label{up1}
\frac{2}{\mu_1^2} g(-h_1)&=&\frac{2}{\mu_1^2} \biggl(\log
\frac{N
\mu
_1^2}{2}+\log\gamma-1+o(1) \biggr).
\end{eqnarray}
As a result, we have the following optimality\vspace*{-2pt} result.
%
%
\begin{thm}\label{thm1}
Assume that the drift sizes are all known and $\mu_i=\mu>0$ for all
$i=1,\ldots,N$.
Then for any predictable, nonsingular, stochastic instantaneous
correlation matrix covariance matrix $\Sigma_t$,
the $N$-CUSUM stopping rule $T_\hbar$ defined in Algorithm \ref
{NCUSUM}, where the set of thresholds
$\hbar$ is chosen using
(\ref{hi}) and~(\ref{equalH}), is asymptotically optimal to the problem
(\ref{eqnproblem}).
More specifically, the difference between the detection delay of the
$N$-CUSUM stopping rule, $J^{(N)}(T_{\hbar})$, and the optimal
detection delay $\inf_{T\in\mathcal{T}_\gamma}J^{(N)}(T)$, is bounded
above by $\frac{2}{\mu_1^2}\log N$, as \mbox{$\gamma\to\infty$}.
\end{thm}

\begin{pf}
The result follows from (\ref{UBJ}), (\ref{lowbound1}) and (\ref{up1}):
\begin{eqnarray*}
0 &\le& J^{(N)}(T_{\hbar}) - \inf_{T\in\mathcal{T}_\gamma}
J^{(N)}(T) \le\frac{2}{\mu_1^2}g(-h_1)-\frac{2}{\mu_1^2}g
\bigl(-\nu_1^\star \bigr)
\nonumber
\\[-2pt]
& \le& \frac{2}{\mu_1^2} \biggl(\log\frac{N \mu_1^2}{2}+\log\gamma- 1 + o(1)
\biggr) - \frac{2}{\mu_1^2} \biggl(\log\frac{\mu_1^2}{2}+\log \gamma -1+o(1)
\biggr)
\\[-2pt]
&=&\frac{2}{\mu_1^2}\log N+o(1),
\end{eqnarray*}
as $\gamma\to\infty$.\vspace*{-4pt}
\end{pf}

On the other hand, in the more general case that the drifts are all
known and
$\mu_1=\cdots=\mu_k<\min_{i>k} \mu_i$, using (\ref{hi}), (\ref{inequalH})
and Result 3 in the Appendix of~\cite{HadjZhanPoor}, we obtain that
%
\begin{eqnarray}
h_1&=&\log\frac{k \mu_1^2}{2}+\log\gamma+o(1),
\\[-2pt]
\label{up2}
\frac{2}{\mu_1^2} g(-h_1)&=&\frac{2}{\mu_1^2} \biggl(\log
\frac{k
\mu
_1^2}{2}+\log\gamma-1+o(1) \biggr).
\end{eqnarray}
It follows that we have the following optimality\vspace*{-2pt} result.\vadjust{\goodbreak}


%
\begin{thm}\label{thm2}
Assume that the drift sizes are known and
$0<\mu
_1=\cdots=\mu_k<\min_{i>k}\mu_i$. Then for any predictable,
nonsingular, stochastic instantaneous correlation matrix $\Sigma_t$,
the $N$-CUSUM stopping rule $T_\hbar$ defined in Algorithm~\ref{NCUSUM}, where the set of thresholds
$\hbar$ is chosen using (\ref{hi}) and (\ref{inequalH}), is
asymptotically optimal to the problem (\ref{eqnproblem}). More
specifically,\vspace*{1pt} the difference between the detection delay of the
$N$-CUSUM stopping rule, $J^{(N)}(T_\hbar)$, and the optimal detection
delay $\inf_{T\in\mathcal{T}_\gamma}J^{(N)}(T)$, is bounded above by
$\frac{2}{\mu_1^2}\log k$, as $\gamma\to\infty$. In particular, if
$k=1$, then $T_\hbar$ is equivalent to the optimal solution to (\ref{eqnproblem})\vspace*{-2pt} asymptotically.
\end{thm}

%
\begin{pf}
The result follows from (\ref{UBJ}), (\ref{lowbound1}) and (\ref{up2}):
\begin{eqnarray*}
0 &\le&J^{(N)}(T_{\hbar}) - \inf_{T\in\mathcal{T}_\gamma}
J^{(N)}(T) \le\frac{2}{\mu_1^2}g(-h_1)-\frac{2}{\mu_1^2}g
\bigl(-\nu_1^\star \bigr)
\nonumber
\\[-2pt]
&\le&\frac{2}{\mu_1^2} \biggl(\log\frac{k \mu_1^2}{2}+\log\gamma- 1 + o(1)
\biggr) - \frac{2}{\mu_1^2} \biggl(\log\frac{\mu_1^2}{2}+\log \gamma -1+o(1)
\biggr)
\\[-2pt]
&=&\frac{2}{\mu_1^2}\log k+o(1),
\end{eqnarray*}
as $\gamma\to\infty$. If $k=1$, the above upper bound for
$J^{(N)}(T_{\hbar}) - \inf_{T\in\mathcal{T}_\gamma} J^{(N)}(T) $ is
$o(1)$, and hence, the $N$-CUSUM stopping rule is equivalent to the
optimal solution to (\ref{eqnproblem})\vspace*{-2pt} asymptotically.
\end{pf}

Finally, if we only know $\mu_1$ and have partial information about
other drifts, that is, $\mu_i\in[\underline{\mu}_i,\overline{\mu}_i]$
for all $i=2,\ldots,N$ and $0<\mu_1=\underline{\mu}_2=\cdots
=\underline{\mu}_{k'}<\min_{i>k'}\underline{\mu}_i$ for some $k'\in\{1,2,\ldots,N-1\}$.
Using (\ref{hi1}), (\ref{inequalH1})
and Result 3 in the Appendix of \cite{HadjZhanPoor}, we obtain that
%
\begin{eqnarray}
h_1& =&\log\frac{\sum_{1\le i\le k'}\mu_1(2\overline
{\mu
}_i-\mu_1)}{2}+\log\gamma+o(1),
\\[-2pt]
\label{up22}
\frac{2}{\mu_1^2} g(-h_1)&=&\frac{2}{\mu_1^2} \biggl(\log
\frac{\sum_{1\le
i\le k'}\mu_1(2\overline{\mu}_i-\mu_1)}{2}+\log\gamma-1+o(1) \biggr).
\end{eqnarray}
It follows that we have the following optimality result.
%
%
\begin{thm}\label{thm3}
Assume that the $\mu_1$ is known, $\mu_i\in
[\underline{\mu}_i,\overline{\mu}_i]$ for all $i=2,\ldots,N$ and that
$0<\mu_1=\underline{\mu}_2=\cdots=\underline{\mu}_{k'}<\min_{i>k'}\underline{\mu}_i$ for some $k'\in\{1,\ldots,\break  N-1\}$. Then for
any predictable, nonsingular, stochastic instantaneous correlation
matrix $\Sigma_t$,
the $N$-CUSUM stopping rule $T_\hbar$ defined in Algorithm~\ref{NCUSUM}, where the set of thresholds
$\hbar$ is chosen using (\ref{hi1}) and (\ref{inequalH1}), is
asymptotically optimal to the problem (\ref{eqnproblem}). More
specifically, the difference between the detection delay of the
$N$-CUSUM stopping rule, $J^{(N)}(T_\hbar)$,\vadjust{\goodbreak} and the optimal detection
delay $\inf_{T\in\mathcal{T}_\gamma}J^{(N)}(T)$, is bounded above by
\[
\frac{2}{\mu_1^2}\log\biggl(\sum_{1\le i\le k'}(2\overline{
\mu}_i/\mu_1-1)\biggr),
\]
as $\gamma\to\infty$. In particular, if $k'=1$, then $T_\hbar$ is
equivalent to the optimal solution to~(\ref{eqnproblem}) asymptotically.
\end{thm}
%
%
\begin{pf}
The result follows from (\ref{UBJ}), (\ref{lowbound1}) and (\ref{up22}):
\begin{eqnarray*}
0 &\le& J^{(N)}(T_{\hbar}) - \inf_{T\in\mathcal{T}_\gamma}
J^{(N)}(T) \le\frac{2}{\mu_1^2}g(-h_1)-\frac{2}{\mu_1^2}g
\bigl(-\nu_1^\star \bigr)
\nonumber
\\
&\le &\frac{2}{\mu_1^2} \biggl(\log\frac{\sum_{1\le i\le k'}\mu
_1(2\overline
{\mu}_i-\mu_1)}{2}+\log\gamma- 1 + o(1)
\biggr) \\
&&{}- \frac{2}{\mu
_1^2} \biggl(\log\frac{\mu_1^2}{2}+\log\gamma-1+o(1)
\biggr)
\\
&=&\frac{2}{\mu_1^2}\log\frac{\sum_{1\le i\le k'}(2\overline{\mu
}_i-\mu
_1)}{\mu_1}+o(1),
\end{eqnarray*}
as $\gamma\to\infty$. If $k'=1$, the above upper bound for
$J^{(N)}(T_{\hbar}) - \inf_{T\in\mathcal{T}_\gamma} J^{(N)}(T) $ is
$o(1)$, and hence, the $N$-CUSUM stopping rule is equivalent to the
optimal solution to (\ref{eqnproblem}) asymptotically.
\end{pf}
%

\begin{rmk}
From Definition \ref{optdef}, we know that the order of the asymptotic
optimality achieved in Theorems \ref{thm1}, \ref{thm2}, \ref{thm3} is
of the second order. If $\mu_1$ is strictly smaller than all the other
drifts ($k=1$ in Theorem \ref{thm2} and $k'=1$ in Theorem \ref{thm3}),
then the $N$-CUSUM stopping rule given by Algorithm \ref{NCUSUM},
and either (\ref{hi}) and (\ref{inequalH}), or (\ref{hi1}) and (\ref{inequalH1})
exhibits third-order asymptotic optimality.
Moreover, it can be seen after a perusal of the proofs that the order
of the asymptotic optimality of the $N$-CUSUM does not change if we
model $\mu_i$'s as $\overline{\mathbb{F}}$-adapted processes bounded by
known constants $\underline{\mu}_i$ and $\overline{\mu}_i$, for all
$i=2,\ldots,N$.
\end{rmk}

\section{Applications}
\label{applications}
In this section, we discuss one of the applications of the results in
decentralized communication systems.
Let us now suppose that each of the observation processes
$\{\xi_t^{(i)}\}$ become sequentially available at a particular
location monitored by
sensor $S_i$, which then employs an asynchronous communication
scheme to a central fusion center. In particular, sensor $S_i$
communicates to the central fusion center only when it wants to
signal an alarm, which is elicited according to a CUSUM stopping rule
$T_{h_i}^{i}$ as in (\ref{CUSUM1chart}) adapted to the small filtration
$\{\mathcal{G}_t^i\}$. The observations
received at the $N$ sensors can change dynamics
at distinct unknown points $\tau_i$. An example of such a case is
described in
\cite{BassBenvGourMeve} where\vadjust{\goodbreak} the motivation suggested arises in the
health-monitoring
of mechanical, civil and aeronautic structures.
The fusion center, whose
objective is to detect the first time when there is a change in at
least one of the sensors devises a minimal strategy; that is, it
declares that a change has
occurred at the first instance when one of the sensors communicates
an alarm. The implication of the main theorems in Section~\ref{Asympt}
is that in fact this strategy is the best, at least
asymptotically, in that there is
no loss in performance, between the case in which the fusion center
receives the raw data $\{\xi_t^{(1)},\ldots\xi_t^{(N)}\}$ summarized
in the large filtration $\{\mathcal{F}_t\}$ directly
and the case in which the communication that takes place is limited to
the decentralized setup.
In other words, the CUSUM stopping rule $T_\hbar$ is a sufficient
statistic (at least asymptotically) of the minimum $N$ possibly
distinct change points. That is, the stopping rule $T_\hbar$ is an
asymptotically optimal solution to the problems of quickest detection
presented in (\ref{eqnproblem}).
In practice sensors are cheap and easy to replace devices whereas
central fusion centers or central processing units are not.
Transferring most of the processing work to the sensors while incurring
no loss in the efficiency of the system is thus
valuable and can render cost and speed effective communication systems.\vspace*{-3pt}

\section{Summary}\label{conclusion}

In this paper, we study the problem of detecting the minimum of $N$
different change points in a $N$-dimensional Brownian system with
partial information of the drifts and an arbitrary, predictable,
nonsingular, stochastic instantaneous correlation matrix $\Sigma_t$. It
is shown that, under an extended Lorden's minmax criterion, the
$N$-CUSUM stopping rule exhibits asymptotic optimality in the tradeoff
between detection delay and false alarms, as the mean time to the first
false alarm increases without bound.
Moreover, the performance of the $N$-CUSUM stopping rule under
dependence is {\em no worse than} that under independence \cite{HadjZhanPoor}.
This optimality result is obtained by establishing a robust upper bound
and a robust lower bound for the optimal detection delay.
The contribution of this work can be seen in two folds. First, we
designed a low complexity, efficient stopping rule without using the
explicit information of the covariance matrix $\Sigma_t$. This stopping
rule is guaranteed to have a comparable performance or identical
performance as the optimal stopping rule, even with cross-correlated
observations---a nontrivial extension to the existing literature and
the first formal treatment of correlated noise in change-point detection.
Second, the robust bounds obtained in this work provide a unified
robust probabilistic (rather than analytical) approach to treat
detection problems with multiple change-points or multiple alternatives
\cite{O.SQA09,HadjMous,ZhanHadjCorr}. This is especially useful when
the analytical characteristic such as joint density or Green functions
are not explicitly available.


\begin{appendix}\label{pf}
\section*{Appendix}

Let us denote by
$\Sigma_t^1=(\rho_t^{i,j})_{i,j\neq1}$ the $(N-1)\times(N-1)$ matrix
obtained from $\Sigma_t$ by removing its first column and first row,\vadjust{\goodbreak}
and by $\tilde{\Sigma}_t^1$ the $(N-1)\times(N-1)$ matrix $(\rho
_{t}^{1,i}\cdot\rho_t^{1,j})_{i,j\neq1}$.
We further introduce a local martingale
\renewcommand{\theequation}{A.\arabic{equation}}
\setcounter{equation}{0}
%
\begin{eqnarray}
\quad B_t^{(1)}&=& -\mu_1 \int
_0^t\bigl(\rho_s^{2,1},
\ldots,\rho_s^{N,1}\bigr) \bigl(\Sigma _s^1-
\tilde {\Sigma}_s^1\bigr)^{-1}
\nonumber
\\[-8pt]
\label{Bt}\\[-8pt]
\nonumber
&&\hspace*{34pt}{}\times \Bigl(\sqrt
{1-\bigl(\rho_s^{2,1}\bigr)^2}\,\dd
\tilde {w}_s^{(2)}, \ldots, \sqrt{1-\bigl(
\rho_s^{N,1}\bigr)^2}\,\dd \tilde
{w}_s^{(N)} \Bigr)',
\\
\label{wtB}
\tilde{w}_t^{(i)}&=& \int
_0^t\frac{\dd w_s^{(i)}-\rho_s^{i,1}\,\dd w_s^{(1)}}{\sqrt{1-(\rho
_s^{i,1})^2}}, \qquad 2\le i\le N.
\end{eqnarray}
The local martingale $B^{(1)}$ will naturally appear in (\ref{eq:last})
of the following proof. We are now ready to prove the assertion in
Lemma \ref{lem:Girsanov}.

\begin{pf*}{Proof of Lemma \ref{lem:Girsanov}}
Since $\Sigma_t$ is nonsingular at all time a.s., we can use a Cholesky
decomposition to obtain a lower triangular, nonsingular matrix-valued
process $L_t=(L_t^{i,j})_{1\le i,j\le N}$, and a $N$-dimensional
standard Brownian motion
$Z=(z^{(1)},\ldots, z^{(N)})'$,
such that
%
\begin{equation}
\dd W_t'=L_t \,\dd Z_t
\end{equation}
holds.
In particular, we have $L_t^{1,1}\equiv1$, $z_t^{(1)}\equiv
w_t^{(1)}$, and $L_t^{i,1}=\rho_t^{i,1}$ for all $2\le  i\le N$.
Using Girsanov's theorem (see, e.g., \cite{KS}, Chapter~3, Theorem~5.1)
and the condition in (\ref{eq:Nov}), the measure changes from $\mathbb
{Q}$ to
$\mathbb{Q}_1$ is given by the exponential martingale
%
\begin{equation}
\frac{d\mathbb{Q}_1}{d\mathbb{Q}} \bigg|_{\mathcal{F}_t}=\mathcal {E}\biggl(\int_0^\cdot
\nu _s\,\dd Z_s\biggr)_t,
\end{equation}
where $\nu=(\nu^{(1)},\ldots,\nu^{(N)})$ is a $N$-tuple process,
such that
%
\begin{equation}
(\mu_1, 0, \ldots, 0)'=L_t\nu'_t.
\end{equation}
It is easily seen that $\nu_t^{(1)}\equiv\mu_1$. Moreover, from
$L_t^{i,1}=\rho_t^{i,1}$ for any $2\le i\le N$, we know that
\begin{equation}
\label{N1}
\tilde{L}_t
\bigl(\nu_t^{(2)},\ldots, \nu_t^{(N)}\bigr)'=-\mu_1 %
\bigl(
\rho_t^{2,1}, \ldots, \rho_t^{N,1} \bigr)',
\end{equation}
%
where $\tilde{L}_t= (L_t^{i,j})_{2\le i,j\le N}$ is a $(N-1)\times
(N-1)$ nonsingular matrix-valued process.
On the other hand, notice that
%
\begin{eqnarray}
&& \dd \bigl(z_t^{(2)},\ldots,
z_t^{(N)}\bigr)'
\nonumber
\\[-8pt]
\label{N2}\\[-8pt]
\nonumber
&& \qquad=(\tilde{L}_t)^{-1}
\bigl(\dd w_t^{(2)}-\rho_{t}^{2,1}\,\dd
w_t^{(1)}, \ldots, \dd w_t^{(N)}-
\rho _{t}^{N,1}\,\dd w_t^{(1)}
\bigr)'.\hspace*{-6pt} 
\end{eqnarray}
%
Using the equations in (\ref{N1}) and (\ref{N2}), we conclude that
%
\begin{eqnarray}
&& \int_0^t\nu_s
\,\dd Z_s
\nonumber\\[-2pt]
&& \qquad =\int_0^t\nu_s^{(1)}
\,\dd z_s^{(1)}+\int_0^t
\bigl(\nu_s^{(2)},\ldots,
\nu_s^{(N)} \bigr)
\bigl(\dd
z_s^{(2)},\ldots, \dd z_s^{(N)}
\bigr)'
\nonumber
\\[-2pt]
&&\qquad =\mu_1 w_t^{(1)}-\mu_1\int
_0^t\bigl(\rho_s^{2,1},
\ldots,\rho _s^{N,1}\bigr) \bigl(\tilde {L}_s'
\bigr)^{-1}(\tilde{L}_s)^{-1}\nonumber\\[-2pt]
&&
\label{eq:last}
\qquad\qquad\hspace*{60pt} {}\times \bigl(\dd
w_s^{(2)}-\rho _{s}^{2,1}\,\dd
w_s^{(1)}, \ldots, \dd w_s^{(N)}-
\rho_{s}^{N,1}\,\dd w_s^{(1)}
\bigr)'
\\[-2pt]
&&\qquad=\mu_1 w_t^{(1)}-\mu_1\int
_0^t\bigl(\rho_s^{2,1},
\ldots,\rho _s^{N,1}\bigr) \bigl(\Sigma _s^1-
\tilde{\Sigma}_s^1\bigr)^{-1}\nonumber\\
&&\hspace*{60pt}\qquad\qquad{}\times \Bigl(\sqrt
{1-\bigl(\rho_s^{2,1}\bigr)^2}\,\dd \tilde {w}_s^{(2)}, \ldots, \sqrt{1-\bigl(
\rho_s^{N,1}\bigr)^2}\,\dd \tilde
{w}_s^{(N)} \Bigr)'
\nonumber
\\[-2pt]
&&\qquad=\mu_1 w_t^{(1)}+B_t^{(1)},
\nonumber
\end{eqnarray}
where the third equality follows from the fact that $\tilde
{L}_t(\tilde
{L}_t)'+\tilde{\Sigma}_t^1=\Sigma_t^1$ holds (see accompanying internet
supplement).
Finally, by the way we construct $B^{(1)}$, we know that the Brownian
motions that drive $B^{(1)}$ and $w^{(1)}$ are independent 
and this completes the proof.
\end{pf*}
\end{appendix}

\section*{Acknowledgments}
The authors are grateful to the Editor and anonymous referees for their
comments contributing to the improvement of this manuscript.






\printaddresses
\end{document}